\documentclass[11pt]{article}
\usepackage{fix-cm}
\usepackage{lmodern}
\usepackage{xr}
\usepackage{amsmath}
\usepackage{amsfonts}
\usepackage{amsthm}
\usepackage{chngcntr}
\usepackage{float}

\usepackage{tikz}
\usetikzlibrary{positioning, shapes.geometric, arrows}
\usepackage{xr}
\usepackage{amsmath}

\makeatletter
\renewenvironment{proof}[1][\proofname]{\par
  \pushQED{\qed}%
  \normalfont \topsep6\p@\@plus6\p@\relax
  \trivlist
  \item[\hskip\labelsep
        \bfseries
        \proofname\ifx#1\proofname\else\space #1\fi. ]}{\popQED\endtrivlist\@endpefalse}
\makeatother

\usepackage{import}

\usepackage{newtxtext}
\usepackage[subscriptcorrection]{newtxmath}

\usepackage{subfiles} 

\graphicspath{{./art/}}

\usepackage[plain,noend]{algorithm2e}

\usepackage{hyperref}
\usepackage{float}

\usepackage{graphics}

\usepackage{newtxtext}
\usepackage[subscriptcorrection]{newtxmath}

\newenvironment{keywords}{\begin{center}\begin{minipage}{0.9\textwidth}\textbf{Keywords:} }{\end{minipage}\end{center}}

\graphicspath{{./art/}}
\usepackage[plain,noend]{algorithm2e}
\makeatletter
\renewcommand{\algocf@captiontext}[2]{#1\algocf@typo. \AlCapFnt{}#2} 
\def\@algocf@capt@plain{top}
\renewcommand{\algocf@makecaption}[2]{%
  \addtolength{\hsize}{\algomargin}%
  \sbox\@tempboxa{\algocf@captiontext{#1}{#2}}%
  \ifdim\wd\@tempboxa >\hsize
    \hskip .5\algomargin%
    \parbox[t]{\hsize}{\algocf@captiontext{#1}{#2}}
  \else%
    \global\@minipagefalse%
    \hbox to\hsize{\box\@tempboxa}
  \fi%
  \addtolength{\hsize}{-\algomargin}%
}
\makeatother

\DeclareMathOperator*{\argmin}{argmin}

\DeclareRobustCommand{\fre}{Fr\'{e}chet }
\newcommand{\cost}{\mathfrak{c}}

\usepackage{adjustbox}

\usepackage{latexsym, amssymb, amscd, amsxtra, amsmath}
\usepackage{graphics, graphicx, color}
\usepackage{natbib}
\usepackage{ifpdf}
\usepackage[format=hang,indention=-1cm,small]{caption}
\usepackage[caption=false]{subfig}
\usepackage{multirow}
\usepackage{hyperref}

\usepackage{apptools}

\newtheorem{thm}{Theorem}
\newtheorem{theorem}[thm]{Theorem}

\theoremstyle{definition}
\newtheorem{definition}{Definition}
\newtheorem*{defn*}{Definition}
\newtheorem{lemma}{Lemma}
\newtheorem{condition}{Condition}
\newtheorem{corollary}{Corollary}
\theoremstyle{remark}
\newtheorem{remark}{Remark}

\AtAppendix{\counterwithin{thm}{section}}  
\AtAppendix{\counterwithin{lemma}{section}}
\AtAppendix{\counterwithin{definition}{section}}

\def\argmin{\mathop{\rm argmin}}

\def\1v{\mathbf 1}
\def\0v{\mathbf 0}

\usepackage{tikz-cd}

\title{Generalized Fr\'{e}chet means with random minimizing domains and its strong consistency}
\author{Jaesung Park and Sungkyu Jung \footnote{The authors were supported by the Samsung Science and Technology Foundation (SSTF-BA2002-03) and the National Research Foundation of Korea (RS-2023-00301976, RS-2024-00333399)} \\ 
Seoul National University}
\begin{document}
\maketitle	

\begin{abstract}
This paper introduces a novel extension of \fre means, referred to as \textit{generalized \fre means}, as a comprehensive framework for describing the characteristics of random elements. The generalized \fre mean is defined as the minimizer of a cost function, and the framework encompasses various extensions of \fre means that have appeared in the literature. The most distinctive feature of the proposed framework is that it allows the domain of minimization for the empirical generalized \fre means to be random and different from that of its population counterpart. This flexibility broadens the applicability of the \fre mean framework to various statistical scenarios, including sequential dimension reduction for non-Euclidean data. We establish a strong consistency theorem for generalized \fre means and demonstrate the utility of the proposed framework by verifying the consistency of principal geodesic analysis on the hypersphere.
\end{abstract}


\section{Introduction}\label{seC:intro}

Many modern statistical tasks can be framed as finding the minimizer of an expected cost (or loss) function. Consider a random element $X$ in a data space $T$ with distribution $P$, and parameter space $M$. The cost function ${\cost}:  T\times M\to\mathbb {R}$ quantifies the cost associated with using a parameter $  m\in M$ to describe a data point $ t\in T$. Often, a key characteristic of the distribution of $X$ is implicitly defined as the parameter that minimizes the expected cost function over a subset $M_{0}$ of $M$. This characteristic is expressed as $\argmin_{m \in M_{0}} E\{ \cost(X, m)\}$, representing the primary target for estimation.

We address the important case where the minimizing domain $M_{0}$ depends on the distribution $P$ and is typically estimated from the data. 
Many statistical estimation problems, such as sequential estimation and constrained minimization, can be framed within this context. A simple example is the estimation of a population tangent vector $m \in M_{0}$, where $M_{0}$ is the tangent space of a Riemannian manifold $T$ at the \fre mean $\mu$ of $P$. Since $\mu$ depends on $P$ and is unknown in general, the domain $M_{0}$ is typically replaced by an empirical domain $M_n$, the tangent space at the sample \fre mean $\hat{\mu}$. 
Another example is sequential dimension reduction in non-Euclidean spaces, where the domain of minimization $M_n$ at each step depends on the results of the previous step, and is dependent to the data and thus random
\citep{fletcher2004principal,jung2012analysis, damon2014backwards, bigot2017geodesic, huckemann2018backward, KIM2024107989}. Since the parameter is a subset of the data space $T$, the parameter space $M$ differs from $T$; see Section~\ref{sec:all_applications} for a specific example. 
Other examples include minimization over a (possibly data-dependent) \emph{restricted} domain such as the $k$-medoids clustering \citep{kaufman2009finding}, which confines the empirical medoids to lie on the empirical support $\{X_1,\ldots, X_n\}$. In this case, the parameter space $M$ is the collection of $k$-tuple subsets of the data space $T$. 

In these situations, the parameter of interest is naturally estimated by minimizing the empirical risk $n^{-1}\sum_{i=1}^{n} \cost(X_{i}, m)$ over a data-dependent, random domain $M_n$. 

This paper introduces a new framework of generalized \fre means with a \emph{random minimizing domain}, allowing the domain  $M_{n}$ of empirical M-estimation to be random. 
Our investigation includes the special case $M_{0} = M_n$, under which the generalized \fre means reduces to a traditional M-estimator. The proposed  framework builds upon and extends the \fre mean \citep{frechet1948elements} and its various extensions: the geometric median \citep{arnaudon2012medians}, the $L^p$ center of mass (for $p \in [1,\infty)$) \citep{afsari2011riemannian}, the $\cost$-\fre mean  \citep{schotz2022strong}, the \fre $\rho$-mean of \cite{huckemann2011intrinsic}, and the domain-restricted \fre means discussed in \cite{evans2024limit}. The generalized \fre mean framework includes these as special cases.
 
The goal of this work is to characterize the general conditions under which a strong law of large numbers, or strong consistency, of generalized \fre mean is guaranteed. We provide a list of minimal conditions essential for consistency, and offer specific, easier-to-check conditions. Our requirements not only ensure consistency for the cases of  random minimizing domain, but are also comparable to, or more relaxed than, previous results in the literature. 

The broad applicability of the framework and the utility of consistency theorems for generalized \fre means are demonstrated by guaranteeing consistency of principal geodesic analysis \citep{fletcher2004principal} for spherical data.

\section{Generalized \fre means}\label{seC:GFM}

\subsection{Definition}
We formally introduce our generalized \fre mean framework. Let $T$ be a topological space where data lie, called \textit{data space}. $X$ denotes a random element from underlying probability space $(\Omega, \mathcal{A}, \mathcal{P})$ to $T$, with distribution $P$. Write $X_{1}, X_{2}, \ldots, X_{n}$ for independent and identically distributed random samples from $P$. We are concerned with certain characteristics of $P$ described by an element in a metric space $(M, d)$ called \textit{descriptor space} or \textit{parameter space}. Note that $M$ may coincide with $T$, but it is generally understood that $M\neq T$. A \textit{cost function} $\cost: T\times M\rightarrow \mathbb{R}$ prescribes the cost $\cost(t,m)$ of using a descriptor $m\in M$ in describing a data point $t\in T$. Throughout this paper, we assume that $\cost(X,m)$ is integrable for all $m\in M$. The population and empirical risk functions, associated with the cost, are  $F(m)=E\{\cost(X,m)\}$ and $F_{n}(m)=\frac{1}{n}\sum_{i=1}^{n} \cost(X_{i},m)$, for $m\in M$, respectively.

Let ${CL}(M)$ be the collection of all nonempty closed subsets of $M$. Let $M_{0} \in {CL}(M)$ be given, possibly depending on the distribution $P$, and let $M_{n}:(\Omega,\mathcal{A},\mathcal{P})\to{CL}(M)$ be a ${CL}(M)$-valued random set for each sample size $n$. 
As mentioned in Section~\ref{seC:intro}, $M_{n}$ can be viewed as an estimator of the true domain $M_{0}$, and may (or may not) depend on the data $\{X_{1}, \ldots, X_{n}\}$.
We defer the discussion on the choice of the $\sigma$-algebra on $CL(M)$ and related measurability issues to the Appendix~\ref{app: sigma algebra}. For any $f: M \rightarrow \mathbb{R}$, $M'\in  CL(M)$ and $\epsilon \ge 0$, we define $\epsilon\text{-}\argmin_{m\in M'} f = \{m\in M' :  f(m) \le \epsilon+\inf_{m\in M'}f(m)\}$, to represent the set of minimizers of $f$ over $M'$, allowing for an $\epsilon$-error.

\begin{definition}\label{def: generalized Frechet mean}
The set of \textit{population generalized \fre means} of $X$ is $E_{0}=\operatorname{argmin}_{m \in M_{0}}~ F(m)$. Similarly, the set of \textit{empirical generalized \fre means} of $X_{1}, X_{2}, \ldots, X_{n}$ is $\hat{E}_{n}=\epsilon_{n}\text{-}{\arg \min }_{m \in M_{n}}~F_{n}(m)$, where $\epsilon_{n} : (\Omega, \mathcal{A}, \mathcal{P}) \rightarrow [0,\infty)$ is a sequence of nonnegative random variables converging to $0$ almost surely. 
\end{definition}

Allowing an $\epsilon_{n}$-error  in the empirical mean accommodates approximate numerical solutions returned by iterative algorithms or optimizers; when an exact minimizer is available, one can simply set $\epsilon_{n}=0$. The sequence $\epsilon_n$ may be either deterministic or data-dependent.

The generalized \fre mean framework is highly flexible, as it does not assume any specific form or structure for the data space $T$, parameter space $M$, and cost function $\cost$. This flexibility enables its application in diverse statistical settings, extending beyond Euclidean spaces to non-Euclidean contexts, such as dimension reduction on manifolds; see Section~\ref{sec:all_applications} and \cite{huckemann2011intrinsic,jung2012analysis,bigot2017geodesic, KIM2024107989}. Additionally, the framework allows the cost function to include penalty terms for $m$ for regularization, while robust loss functions such as the Huber loss \citep{Huber1964} can be used to improve robustness.



\subsection{Relations to previous extensions of \fre means}\label{seC:frechet means}

Our framework encompasses various extensions of \fre means as special cases.

For $T \neq M$, where the parameter $m$ lies in the data space,  the \fre mean and the $L^p$ center of mass~\citep{afsari2011riemannian} serve as location parameters for the distribution $P$. The population and empirical $L^p$ centers of mass are given by $\argmin_{m\in M}~ E \{d^{p}(X,m)\}$ and ${  \operatorname{argmin}}_{m\in M}~\frac{1}{n}\sum_{i=1}^{n} d^{p}(X_{i},m)$,  respectively. Special cases include the geometric median ($p=1$) and \fre mean ($p=2$). Recently, \cite{evans2024limit} introduced 
$C$-restricted \fre $p$-mean, a generalization of the $L^{p}$ center of mass that restricts the minimizing domain $C$ to $M_{0}$ for the population case and random $M_n$ in the empirical case.  

The case $T \neq M$ was first addressed by  \cite{huckemann2011intrinsic}, who generalized the squared loss $d^{2}(X,m)$ of the \fre mean to $\rho^{2}(X,m)$, where $\rho: T\times M \rightarrow [0,\infty)$ acts as  a ``distance'' function between a data point $X$ and a descriptor $m$ that lie in different spaces. More recently, \cite{schotz2022strong} extended this by replacing $\rho^{2}$ with a more flexible class of cost functions $\cost$, which can take negative values. Both \fre $\rho$-means of \cite{huckemann2011intrinsic} and $\cost$-\fre means of \cite{schotz2022strong}  are also applicable when $T= M$. In a similar setting, \cite{brunel2023geodesically} considered a geodesically convex cost function $\cost$. 


Table~\ref{Table: Frechet means} summarizes the extent to which \fre means have been generalized. The proposed generalized \fre mean framework encompasses all aspects of generalization, including different data space and parameter space, flexible cost functions, and random minimizing domains, all of which are essential for the application demonstrated in Section \ref{sec:all_applications}. When $M_{0}=M_{n}=M$, the empirical generalized \fre mean coincides with the general empirical M-estimator.
 



\begin{table}[!h]
\caption{Extensions of \fre means}
\label{Table: Frechet means}
\resizebox{\textwidth}{!}{%
\begin{tabular}{cccc}
                              & Data space $T$           & Objective function             & Minimizing domain \\
\fre mean~\citep{frechet1948elements}             & $T=M$                    & $d^{2}(t,m)$                    & $M_{0}=M_{n}=M$                  \\
$L^{p}$ center of mass~\citep{afsari2011riemannian}        & $T=M$                    & $d^{p}(t,m), (p\ge 1)$                    & $M_{0}=M_{n}=M$                  \\
\fre $\rho$-mean~\citep{huckemann2011intrinsic}       & Possibly $T\neq M$    & $\rho^{2}(t,m)$                 & $M_{0}=M_{n}=M$                  \\
$\cost$-\fre mean~\citep{schotz2022strong}      & Possibly $T\neq M$    & $\cost(t,m)\rightarrow \mathbb{R}$                    & $M_{0}=M_{n}=M$                  \\
restricted \fre $p$-mean~\citep{evans2024limit} & $T=M$ & $d^{p}(t,m), (p\ge 1)$ &  Possibly $M_{0} \neq M_{n} \subset M$ \\
Generalized \fre mean  & Possibly $T\neq M$    & $\cost(t,m)\rightarrow \mathbb{R}$          & Possibly $M_{0} \neq M_{n} \subset M$  \\
\end{tabular}
}
\end{table}

\subsection{Strong consistency of generalized \fre means}\label{sec:theory}
This subsection introduces sufficient conditions for ensuring strong consistency of the empirical generalized \fre mean set, where the minimization is performed over potentially random domains $M_n$. Intuitively, this requires some form of convergence of $M_n$ to $M_0$. Recall that ${CL}(M)$ denotes the collection of all non-empty closed subsets of $M$.

\begin{definition}\label{Kuratowski}\citep{beer1993topologies}
For \(B_{0} \in {CL}(M)\) and \(B_{n} \in {CL}(M)\), we say that \(B_{n}\) converges to \(B_{0}\) in the sense of Kuratowski if the following conditions are satisfied. 1) Any accumulation point of an arbitrary sequence \(m_{n} \in B_{n}\) lies in $B_0$. 2) For every \(m_{0} \in B_{0}\), there exists a sequence \(m_{n} \in B_{n}\) that converges to \(m_{0}\).
\end{definition}

\begin{condition}\label{C:Kuratowski}
$M_n$ converges to \(M_{0}\) in the sense of Kuratowski almost surely.
\end{condition}

We also give conditions for the parameter space $M$ as follows.

\begin{condition}\label{C:separable}
$(M,d)$ is separable and complete. In other words, $M$ is Polish.
\end{condition}

\begin{condition}\label{C:Heine-Borel}
Every closed and bounded subset of $(M,d)$ is compact. In other words, $M$ satisfies the Heine–Borel property.
\end{condition}

In addition, we specify the conditions required for the cost function \(\cost\). While prior works such as \cite{schotz2022strong, jung2025averaging}, which consider fixed domains, rely only on lower semi-continuity of \(\cost(t,\cdot)\), upper semi-continuity becomes essential when the minimizing domain itself is random, as in our framework.

\begin{condition}\label{C:rho-continuous}
For each $t \in T$, $\cost(t, \cdot): M \rightarrow \mathbb{R}$ is continuous.
\end{condition}

\begin{definition}\label{def: local-inf-sup}
For $m\in M$ and $r>0$, we denote local infimum and local supremum of cost function, $\pi_{m,r}: T\rightarrow \mathbb{R}$ and $\Pi_{m,r}: T \rightarrow \mathbb{R}$ by $\pi_{m,r}(t)= \inf_{m' \in B(m, r)} \cost(t, m')$, $\Pi_{m,r}(t)= \sup_{m' \in B(m, r)} \cost(t, m')$, where $B(m, r) \subset M$ denotes the open ball with center $m$ and radius $r > 0$.  
\end{definition}

\begin{condition}\label{C:locally-integrable}
For every $m \in M$, there exists $r_{m} > 0$ such that $\pi_{m,r}(X)$ and $\Pi_{m,r}(X)$ are integrable for all $0<r\le r_{m}$. 
\end{condition}

\begin{theorem}\label{thm: BP SC}
Suppose that Conditions~\ref{C:Kuratowski}--\ref{C:locally-integrable} hold, and that $\hat{E}_{n}$ is \emph{almost surely eventually bounded}; that is, with probability one, there exists \(N \ge 1\) such that \(\cup_{n \ge N} \hat{E}_{n}\) is bounded. Then, 
\begin{equation}\label{eq: BP consistency}
\mathcal{P} \left[ \left\{\omega\in  \Omega  :   \lim_{n\rightarrow \infty}  \textstyle \sup_{m\in \hat{E}_{n}(\omega)} d(m, E_{0}) =0 \right\} \right]=1.    
\end{equation}
\end{theorem}
For the proof of Theorem~\ref{thm: BP SC}, see Appendix~\ref{app: technical details for sec 2.3}. We say that $\hat{E}_{n}$ is a strongly consistent estimator of $E_{0}$ in the sense of Bhattacharya-Patrangenaru (or $\hat{E}_{n}$ is BP consistent with $E_{0}$) when \eqref{eq: BP consistency} holds \citep{Bhattacharya2003}. When $E_{0}$ is evident from the context, we simply refer to $\hat{E}_n$ as BP consistent.

BP consistency ensures that every selection $e_n\in\hat E_n$ converges almost surely to $E_0$, but it does not guarantee that for each $e\in E_0$  there exists a sequence $e_{n}\in \hat{E}_{n}$ converging to that specific $e$; the convergence is thus often referred as one-sided \citep{schotz2022strong}. Stronger convergence in the Hausdorff distance, $h_d(\hat E_n,E_0)=\max\{\sup_{m\in\hat E_n} d(m,E_0),\,\sup_{m\in E_0} d(m,\hat E_n)\}$, is more desirable but typically harder to verify. See \cite{blanchard2025frechet} for a result on convergence in Hausdorff distance under a special case $T=M$. Nevertheless, under the natural assumption that $E_{0} = \{e_{0}\}$ is a singleton, BP consistency guarantees convergence in Hausdorff distance.

\begin{lemma}\label{lem: singletone}
Let $\hat{E}_{n}$ be BP consistent with $E_{0}$, and suppose that $E_{0}$ consists only of a single element $e_{0}$. Then, almost surely, all sequences $e_{n}\in \hat{E}_{n}$ converge to $e_{0}$.
\end{lemma}

\begin{proof}[of Lemma~\ref{lem: singletone}]
Since $E_{n}$ is BP consistent with $E_{0}=\{e_{0}\}$, we observe that for almost all $\omega$, ${\lim}_{n\rightarrow \infty} {\sup}_{e\in {E}_{n}(\omega)} d(e, e_{0}) =0$ holds. Thus, ${\lim}_{n\rightarrow \infty}d(e_{n}, e_{0}) =0$ for all sequence $e_{n}\in {E}_{n}(\omega)$ and for almost all $\omega$.
\end{proof}

We now briefly explain Conditions~\ref{C:Kuratowski}--\ref{C:locally-integrable}, as well as the notion of almost sure eventual boundedness of \(\hat{E}_{n}\), and discuss how these assumptions ensure BP consistency of $\hat{E}_{n}$.

Condition~\ref{C:Kuratowski}, requiring $M_n\to M_0$ in the Kuratowski sense, holds trivially when $M_n\equiv M_0$, but must be verified separately in applications involving genuinely random domains, as discussed in Section~\ref{sec:all_applications}. Kuratowski convergence is weaker than Hausdorff convergence: the latter implies the former, but not conversely \citep[Remark~1.3]{schotz2022strong}. The two notions coincide when $M$ is compact \citep[Chapter~5]{beer1993topologies}. When $M$ is unbounded, $h_{d}(M_n,M_0)$ may be infinite for all $n$, so requiring Hausdorff convergence is often impractical.


Conditions~\ref{C:separable} and \ref{C:Heine-Borel} are both standard in the analysis of consistency for empirical minimizers \citep{ziezold1977expected, korf2001random, afsari2011riemannian, huckemann2011intrinsic, schotz2022strong, evans2024limit}. While completeness is not strictly required in some earlier works \citep{frechet1948elements, huckemann2011intrinsic, evans2024limit}, it is essential in our setting to ensure the measurability of the local infimum and supremum functions $\pi_{m,r}$ and $\Pi_{m,r}$ introduced in Definition~\ref{def: local-inf-sup}; see \cite{korf2001random} for further details. Since every metric space admits a completion \citep[Theorem 1.15]{howes2012modern}, it is natural to assume completeness of $M$.

Conditions~\ref{C:rho-continuous} and \ref{C:locally-integrable} are technical assumptions that, when combined with Conditions~\ref{C:Kuratowski} and \ref{C:separable}, imply the following:
\begin{equation}\label{eq: Z-SC}
\mathcal{P} \left[ \left\{\omega\in  \Omega  : \cap_{n=1}^{\infty} \overline{\cup_{k=n}^{\infty} \hat{E}_{n}(\omega)} \subset E_{0}  \right\} \right]=1.
\end{equation} 
When $\hat{E}_{n}$ satisfies \eqref{eq: Z-SC},  $\hat{E}_{n}$ is said to be Ziezold consistent with $E_{0}$ \citep{ziezold1977expected}. Ziezold consistency implies that almost surely any accumulation point of a sequence chosen from $\hat{E}_n$ lies within $E_{0}$. While Ziezold consistency is strictly weaker than BP consistency, it implies BP consistency when \(\hat{E}_n\) contains no diverging subsequence, as discussed in \cite{huckemann2011intrinsic, schotz2022strong}. The almost sure eventual boundedness of \(\hat{E}_n\) assumed in  Theorem~\ref{thm: BP SC}, together with Condition~\ref{C:Heine-Borel}, rules out diverging subsequences and implies BP consistency.

\begin{remark}\label{rmk:transfer-ergodic}
A key technical step in verifying Theorem~\ref{thm: BP SC} is to show $F_{n}(m_{n}) \rightarrow F(m_{0})$, which is essential in establishing Ziezold consistency \eqref{eq: Z-SC}, for $m_{n}\in M_{n}$ converging to $m_{0}\in M_{0}$, \emph{without} assuming a special form or strong regularity of the cost function. 
In our proof, this convergence is achieved through the convergence of the averages of local envelopes $\pi_{m,r}$ and $\Pi_{m,r}$, which allows us to accommodate a broad family of cost functions. Moreover, we believe that our envelope-based approach can be extended to verify convergence of generalized \fre means from non-i.i.d. data.
For example, we conjecture that Birkhoff’s ergodic theorem \citep[Sec~1.6]{walters2000introduction} provides the required limits for the relevant envelope averages for ergodic sequences under certain conditions. A systematic treatment 
is left for future work.



\end{remark}

In the next subsection, we present examples of cost functions that satisfy Conditions~\ref{C:rho-continuous} and \ref{C:locally-integrable}. The almost sure eventual boundedness of $\hat{E}_n$ is also established for specific classes of cost functions and distributions, building on the frameworks developed in \cite{huckemann2011intrinsic} and \cite{schotz2022strong}.

\subsection{Examples of cost functions satisfying BP consistency}\label{seC:various examples}
We highlight specific classes of cost functions that satisfy Conditions~\ref{C:rho-continuous} and \ref{C:locally-integrable}, and ensure
the eventual boundedness of $\hat{E}_n$. These classes include previously proposed extensions of \fre means, as discussed 
in 
\cite{frechet1948elements,afsari2011riemannian,huckemann2011intrinsic,evans2024limit}.

We first consider the general case, where $M$ and $T$ may differ. The class of cost functions we introduce extends those used in \cite{huckemann2011intrinsic}, where the cost function takes the form $\cost(t,m) = \rho^2(t,m)$. 
In our framework, the cost functions are generalized to $\cost(t,m) = G\{ \rho(t,m)\}$, where $\rho$ and $G$ satisfy certain regularity conditions.
 
\begin{definition}[\cite{huckemann2011intrinsic}]\label{def: equicontinuous} \label{def: coercive} 
For $\rho: T \times M \to [0,\infty)$, we define: 

(i) \emph{Equicontinuity}: $\rho$ is equicontinuous if, for every $m \in M$ and $\epsilon>0$, there exists $\delta(\epsilon, m)>0$ such that $|\rho(t, m^{\prime})-\rho(t, m)|<\epsilon$ for all $t\in T$ and $m^{\prime} \in M$ with $d( m, m^{\prime})<\delta(\epsilon, m)$. 

(ii) \emph{Coercivity}: $\rho$ is coercive (with respect to $P$) if there exist $o \in M$ and $C > 0$ such that:  

1. $\mathcal{P}\{ \rho\left(X, o \right) <C \}>0$, 

2. For any sequence $m_{n} \in M$ with $d(o, m_{n}) \to \infty$, there exists a real sequence $A_{n} \to \infty$ such that $\rho(t, m_{n})>A_{n}$ for all $t \in T$ satisfying $\rho(t, o)<C$.
\end{definition}
 


\begin{corollary}\label{cor: BP SC huckemann generalization}
Suppose that the cost function satisfies $\cost=G\circ \rho$, where $G: [0,\infty) \rightarrow [0,\infty)$ is non-decreasing, continuous, \emph{$b$-subadditive} for some $b\ge 1$ and satisfies $\lim_{x\rightarrow \infty} G(x)=\infty$, and $\rho: T\times M \rightarrow [0,\infty)$ is equicontinuous and coercive. If Conditions~\ref{C:Kuratowski}, \ref{C:separable}, and \ref{C:Heine-Borel} hold, then \(\hat{E}_{n}\) is BP consistent with \(E_{0}\).
\end{corollary}

For the proof of Corollary~\ref{cor: BP SC huckemann generalization}, see Appendix~\ref{app: technical details for Sec 2.4}.

In Corollary~\ref{cor: BP SC huckemann generalization}, a nondecreasing function $G$ is called \emph{$b$-subadditive} if $G(2x)\le bG(x)$ for all $x\ge0$. Concave and non-decreasing functions or polynomials of the form $G(x)=x^{p}$ ($p\ge1$) satisfy $b$-subadditivity for some $b\ge 1$, while exponentially growing functions do not. 
The subadditivity and equicontinuity controls 
the local behavior of $\cost(t,m')$ for $m'$ near $m$, and ensures Conditions~\ref{C:rho-continuous} and \ref{C:locally-integrable}. On the other hand, the condition $\lim_{x\to\infty}G(x)=\infty$ together with the coercivity 
implies that $F_n(m)\to\infty$ as $d(m,o)\to\infty$, guaranteeing the 
eventual boundedness of $\hat E_n$.

 

\begin{remark}\label{rmk: equicont cost function} 

Equicontinuous cost functions frequently arise in statistical problems involving projections in metric spaces. A key example is dimension reduction for data in a Riemannian manifold $(T,d_T)$, where the goal is to find the $k$-dimensional submanifold $m \subset T$ that minimizes the sum of squared projection errors. Here, the parameter space $(M,d)$ is the collection of nonempty subsets of $T$, endowed with the Hausdorff distance $d = h_{d_T}$ induced by $d_T$. The typical cost function $\cost(t,m) = d_T^2(t,m)$, where $d_T(t,m) = \inf_{t' \in m} d_T(t,t')$ is defined as the distance from the point $t \in T$ to the set $m \subset T$, utilizes the metric projection of $t$ onto $m$. 
%
%
For all \(t \in T\) and \(m, m' \in M\), the inequality \(\mid d_{T}(t, m) - d_{T}(t, m') \mid \leq h_{d_T}(m, m')\) holds, establishing that $\rho = d_T$ is equicontinuous. The application in Section~\ref{sec:all_applications} utilizes this notion of metric projection.
\end{remark}

When parameter space $M$ and data space $T$ are equal and $\rho(t,m)=d(t,m)$, $\rho$ trivially satisfies equicontinuity and coercivity in Definition~\ref{def: coercive}, by the triangle inequality of the distance function $d$. We obtain the following result as a direct application of Corollary~\ref{cor: BP SC huckemann generalization}. See Appendix~\ref{app: technical details for Sec 2.4} for the detailed proof of Corollary~\ref{cor: BP SC Hd}.

\begin{corollary}\label{cor: BP SC Hd}
Suppose $T=M$, $\cost=G\circ d$, where $G: [0,\infty) \rightarrow [0,\infty)$ is non-decreasing, $b$-subadditive for some $b\ge 1$, continuous and satisfies $\lim_{x\rightarrow \infty} G(x)=\infty$. If Conditions~\ref{C:Kuratowski}, \ref{C:separable}, and \ref{C:Heine-Borel} hold, then \(\hat{E}_{n}\) is BP consistent with \(E_{0}\). In particular, this conclusion holds for $\cost = d^p$ ($p\ge 1$).
\end{corollary}

As demonstrated at the end of Corollary~\ref{cor: BP SC Hd}, BP consistency of the empirical \(L^p\)-center of mass can be established even when the minimization is performed over a random domain. This result has also been proven in Theorem 4.4 of \cite{evans2024limit}. Therefore, Corollary~\ref{cor: BP SC Hd} can be viewed as a generalization of the corresponding result in \cite{evans2024limit}.

\begin{remark}\label{rmk: centered Lp mean}
In the setting where \(M = T\) and the domain is fixed as \(M_n = M_0 = M\), \cite{schotz2022strong} showed BP consistency of $\hat{E}_{n}$ with a shifted cost function \(\cost(t, m) = d^p(t, m) - d^p(t, o)\), for an arbitrarily chosen reference point \(o \in M\). In the Appendix~\ref{app: further generalization}, we extend this result by establishing BP consistency under the same shifted cost function, but in a more general setting where the minimizing domains $M_n$ may be random. Using this modified cost function allows a weaker moment condition, namely \(E\{d^{p-1}(X, m)\} < \infty\), instead of the usual \(E\{d^p(X, m)\} < \infty\), which has been assumed throughout. Moreover, this consistency result is extended to the case where $0< p \le 1$, without requiring any integrability condition. 
\end{remark}

\section{Application: Principal Geodesic analysis in Hypersphere}\label{sec:all_applications}

The principal geodesic analysis of \cite{fletcher2004principal} is a widely used dimension reduction method for directional data on hyperspheres $\mathbb{S}^{p-1}= \{ x \in \mathbb{R}^{p} : \|x\|_{2} = 1 \}$. Given data $X_1,\ldots,X_n \in \mathbb{S}^{p-1}$, the method sequentially identifies the best-fitting geodesic submanifolds (subspheres) of increasing dimensions. Briefly, the method works as follows: Using the standard \fre mean 
$\hat\mu$ of the data as the base point, the first principal geodesic $\hat{S}^1$ is found among the geodesics passing through $\hat\mu$. Given $\hat{S}^1$, the best-fitting 2-dimensional geodesic submanifold $\hat{S}^2$ is determined among those containing $\hat{S}^1$. This process continues, producing a nested sequence of principal geodesic submanifolds 
\begin{equation}
    \label{eq:samplePGAseq}
    \{\hat\mu\} = \hat{S}^0 \subset \hat{S}^1 \subset \hat{S}^2 \subset \cdots \subset \hat{S}^{p-1} = \mathbb{S}^{p-1}.
\end{equation} 

The estimands in the principal geodesic analysis are defined as follows: For $k = 1,\ldots, p-1$, 
\begin{equation}
    \label{eq:popPGA}
    S^{k} = \argmin_{S \in \mathfrak{S}^{k}(S^{k-1})} E \{d^{2}(X, S)\},
\end{equation} 
where $d(s,S)$ is the geodesic distance between a point $s\in \mathbb{S}^{p-1}$ and a set $ S\subset \mathbb{S}^{p-1}$, 
and $\mathfrak{S}^{k}(S^{k-1})$ is the collection of all $k$-dimensional subspheres of $\mathbb{S}^{p-1}$ that contain $S^{k-1}$. In contrast, the sample-based estimators in  \eqref{eq:samplePGAseq} are given by
$$\hat{S}^{k} = \argmin_{S \in \mathfrak{S}^{k}(\hat{S}^{k-1})} \frac{1}{n} \sum_{i=1}^n d^{2}(X_i, S),$$
where the domain of minimization $\mathfrak{S}^{k}(\hat{S}^{k-1})$ differs from that of \eqref{eq:popPGA} and depends on the estimator $\hat{S}^{k-1}$ in the previous step. 

Our generalized \fre mean framework can be used to establish the consistency of the empirical sequence \eqref{eq:samplePGAseq}, derived from a random sample, with its population counterpart \(\mu \in S^1 \subset \cdots \subset S^{p-1} = \mathbb{S}^{p-1}\). To this end, we identify \(S^k\) and \(\hat{S}^k\) as the population and empirical generalized \fre means, respectively, where the minimizing domains are given by \(M_0 = \mathfrak{S}^k(S^{k-1})\) and \(M_n = \mathfrak{S}^k(\hat{S}^{k-1})\). We obtain the following result by repeatedly applying Corollary~\ref{cor: BP SC Hd} in an inductive manner over increasing \(k\). See Appendix~\ref{app: technical details for Sec 2.5} for the detailed proof of Theorem~\ref{thm: BP consistency of PGA}.

\begin{theorem}\label{thm: BP consistency of PGA}
Under the above setting, suppose the sequence \(\mu \in S^{1} \subset \cdots \subset S^{p-1}\) of principal geodesic subspheres is uniquely defined. Then, as \(n \to \infty\), \(h_d(\hat{S}^k, S^k) \to 0\) almost surely for all \(k = 0, 1, \ldots, p-1\).
\end{theorem}

In Theorem~\ref{thm: BP consistency of PGA}, \(h_d\) denotes the Hausdorff distance induced by \(d\).  To apply Corollary~\ref{cor: BP SC Hd}, we verify Conditions~\ref{C:Kuratowski}--\ref{C:Heine-Borel}. Condition~\ref{C:Kuratowski} holds since the almost sure convergence of \(\hat{S}^{k-1}\) to \(S^{k-1}\) implies the Kuratowski convergence of \(\mathfrak{S}^k(\hat{S}^{k-1})\) to \(\mathfrak{S}^k(S^{k-1})\). The compactness of the sample space \(\mathbb{S}^{p-1}\) and parameter space \(\mathfrak{S}^k\) ensures Conditions~\ref{C:separable} and \ref{C:Heine-Borel}. As noted in Remark~\ref{rmk: equicont cost function}, the cost function based on metric projection satisfies the equicontinuity of $\rho$. The coercivity of $\rho$ is trivial since \(\mathfrak{S}^k\) is compact. This consistency result also extends to cases where \(d^2(x, S)\) is replaced by \(G\{d(x, S)\}\) for any $G$ satisfying the conditions in Corollary~\ref{cor: BP SC huckemann generalization}.

The uniqueness assumption in Theorem~\ref{thm: BP consistency of PGA} is imposed solely for clarity of presentation. Even when multiple population sequences exist, the empirical sequence \(\hat{S}^k\) is still BP consistent, in the sense that it converges almost surely to one of the true minimizers.

\appendix

\section{Technical details and additional results}

This section contains technical details and additional results not included in the main paper.

\subsection{\texorpdfstring{$\sigma$}{sigma}-algebra on \texorpdfstring{$CL(M)$}{CL(M)} and measurability issue for \texorpdfstring{$\hat{E}_{n}$}{En}}\label{app: sigma algebra}

In this subsection, we discuss the choice of the $\sigma$-algebra on ${CL}(M)$ and the measurability issues of $M_{n}$ and $\hat{E}_n$.

As mentioned in \cite{evans2024limit}, a common choice of $\sigma$-algebra on $CL(M)$ is the Effros $\sigma$-algebra \citep{srivastava1998course}, defined by the smallest $\sigma$-algebra containing $$\{ \{C \in CL(M) : C\cap U \neq \emptyset \} : U\subset M \text{ is open } \},$$ and denoted by $\mathcal{E}(X)$. In our development, we assume that $M_{n}(\omega)$ is measurable with respect to the Effros $\sigma$-algebra, but such an assumption can be made for any other choice of the $\sigma$-algebra on $CL(M)$.

It can be easily shown that under Condition $\ref{C:rho-continuous}$, $\hat{E}_{n} \in CL(M)$ for all $\omega\in \Omega$ since $\hat{E}_{n}$ is the inverse image of a closed set of a continuous function. However, it is unclear whether our $\hat{E}_{n} : (\Omega, \mathcal{A}, \mathcal{P}) \rightarrow (CL(M), \mathcal{E}(X))$ is measurable under certain $T$, $M$, $\cost$ and $M_{n}$. This issue is addressed now.


Note that, 
regardless the choice of $\sigma$-algebra on $CL(M)$, our primary event $$ \left\{\omega\in  \Omega  :   \lim_{n\rightarrow \infty}  \textstyle \sup_{m\in \hat{E}_{n}(\omega)} d(m, E_{0}) =0 \right\} $$ 
contains an event of probability one, by Theorem~\ref{thm: BP SC} under the stated conditions. 
Hence, it is measurable for the $\mathcal{P}$-completion \citep[Section A.2]{durrett2019probability} of the $\mathcal{A}$, and the particular choice of $\sigma$-algebra on $CL(M)$ is immaterial in our case.

For a detailed treatment of the measurability issues for $\hat{E}_{n}$ in the special setting $T=M$, $\cost= d^{p}$, with respect to the Effros $\sigma$-algebra, we rely on the development in Section 3 of \cite{evans2024limit}. To discribe the measurability of the sample \fre mean sets more generally, we restate the results of \cite{evans2024limit} using our notation.


\begin{theorem}[Lemma 3.2 of \cite{evans2024limit}]\label{thm: Enhat measurable}
Let $(M,d)$ be a metric space and $T$ is equal to $M$. Let $M_{0} = \text{supp}(P)$ and $M_{n} = \{ X_{1}, X_{2}, \ldots, X_{n}\}$. For $1\le p<\infty$, define $\cost(t,m)= d^{p}(t,m)$ for $t\in T$ and $m\in M$. Then, $\hat{E}_{n} : (\Omega, \mathcal{A}, \mathcal{P}) \rightarrow (CL(M), \mathcal{E}(X))$ is measurable.
\end{theorem}

\begin{theorem}[Lemma 3.3 of \cite{evans2024limit}]\label{thm: Z-SC measurable}
Assume $(X,d)$ is separable and locally compact. For any sequence of measurable maps $C_{n} : (\Omega, \mathcal{A}, \mathcal{P}) \rightarrow (CL(M), \mathcal{E}(M))$, a random closed set $$\cap_{n=1}^{\infty} \overline{\cup_{m=n}^{\infty} C_{m}} : (\Omega, \mathcal{A}, \mathcal{P})  \rightarrow (CL(M), \mathcal{E}(M))$$ is measurable.
\end{theorem}

\begin{theorem}[Lemma 3.4 of \cite{evans2024limit}]\label{thm: one-sided distance measurable}
Let $K(M) \subseteq CL(M)$ be the collection of all nonempty closed compact subsets of $CL(M)$. If $C:(\Omega, \mathcal{A}, \mathcal{P})\rightarrow (K(M), \mathcal{E}(M))$ is measurable and $C'\in CL(M)$ is fixed, then $$\max_{c\in C}\min_{c'\in C'} d(c,c'): (K(M), \mathcal{E}(M), \mathcal{P}) \rightarrow (\mathbb{R}, \mathcal{B}(\mathbb{R}))$$ is measurable, where $\mathcal{B}(\mathbb{R})$ is the standard Borel $\sigma$-algebra on $\mathbb{R}$.  
\end{theorem}

Theorem~\ref{thm: Enhat measurable} guarantees the measurability of \(\hat{E}_{n}\) for specific choices of \(M_{0}\) and \(M_{n}\). When combined with Theorem \ref{thm: Z-SC measurable}, it ensures that the event \[ \left\{w \in \Omega : \cap_{n=1}^{\infty} \overline{\cup_{m=n}^{\infty} \hat{E}_{m}} \subset E_{0} \right\}\] is measurable. In other words, it guarantees that the event of \(\hat{E}_{n}\) being Ziezold consistent with respect to \(E_{0}\) is measurable without the need for completion.

Finally, by combining Theorem~\ref{thm: one-sided distance measurable} and Theorem~\ref{thm: Enhat measurable}, it is ensured that the set \[ \left\{\omega \in \Omega : \lim_{n \to \infty} \max_{e \in \hat{E}_{n}(\omega)} \min_{e'\in E_{0}} d(e, e') = 0\right\}\] is measurable, when $\hat{E}_{n}$and $E_{0}$ are all compact. This can be guaranteed when $(M,d)$ satisfies the Heine-Borel property, since $F_{n}(m)$ and $F(m)$ diverges as $d(m,o)\rightarrow \infty$ for any fixed $m$, thus $\hat{E}_{n}$ and $E_{0}$ are closed, bounded, and compact under Heine-Borel property. This guarantees that the event where \(\hat{E}_{n}\) is BP consistent with respect to \(E_{0}\) is also measurable without the need for completion.

\subsection{Technical details for Section~\ref{sec:theory}}\label{app: technical details for sec 2.3}

In this section, we establish Theorem~\ref{thm: BP SC} and introduce several technical lemmas that support its proof. To clarify the logical flow, we outline the structure and progression of our main technical result in Section~\ref{sec:theory}, illustrated by the diagram in Figure~\ref{fig:consistency_diagram}.

We prove BP consistency (Lemma~\ref{lem:Z-SC-BP-SC-WC}) by reducing it to Ziezold consistency and eventual boundedness (Condition~\ref{C:Heine-Borel}); Ziezold consistency follows from a.s. continuous convergence (Lemma~\ref{lem:epi-convergence-to-Z-SC}), which is delivered by the Lemma~\ref{lem:epi-convergence}) under Conditions~\ref{C:separable}, \ref{C:rho-continuous}, and \ref{C:locally-integrable}. Then we finally get Theorem~\ref{thm: BP SC}, BP consistency of $\hat{E}_{n}$ under Conditions~\ref{C:Kuratowski}--\ref{C:locally-integrable} and almost sure eventual boundedness of $\hat{E}_{n}$.

\begin{figure}
    \centering
    \begin{tikzpicture}[
        node distance=1cm and 0.5cm,
        every node/.style={
            draw,
            text width=2.9cm,
            align=center,
            rounded corners,
            font=\small,
            fill=gray!10
        },
        every path/.style={->, thick}
    ]
        \node (BP)              [rectangle] at (0, 0)   {BP Consistency \\ \textit{(Lemma~\ref{lem:Z-SC-BP-SC-WC})}};
        \node (Ziezold)         [rectangle] at (0,-2)   {Ziezold Consistency \\ \textit{(Lemma~\ref{lem:epi-convergence-to-Z-SC} )}};
        \node (FMConverge)      [rectangle] at (0,-4)   {$F_n$ continuously converges to $F$ almost surely \\ \textit{(Lemma \ref{lem:epi-convergence})}};
        \node (Conditions245)   [rectangle] at (0,-6)   {Conditions \ref{C:separable}, \ref{C:rho-continuous}, \ref{C:locally-integrable}};
        \node (Condition3)      [rectangle] at (4,-2)   {Condition~\ref{C:Heine-Borel}};
        \node (EventuallyBounded)[rectangle] at (8,-2)  {$\hat{E}_n$ is eventually bounded almost surely};
        \node (Condition1)      [rectangle] at (-4,-4)  {Condition \ref{C:Kuratowski}};

        \draw (Ziezold) -- (BP);
        \draw (Condition3) -- (BP);
        \draw (EventuallyBounded) -- (BP);
        \draw (FMConverge) -- (Ziezold);
        \draw (Condition1) -- (Ziezold);
        \draw (Conditions245) -- (FMConverge);
    \end{tikzpicture}
    \caption{Conditions ensuring BP consistency}
    \label{fig:consistency_diagram}
\end{figure}

As mentioned in Section~\ref{sec:theory}, $\hat{E}_{n}$ is BP consistent with $E_{0}$ when it is Ziezold consistent and there are no diverging sequences. We formally state this result as the following Lemma.

\begin{lemma}[Theorem 2.4 of \cite{schotz2022strong}]\label{lem:Z-SC-BP-SC-WC}
If $\hat{E}_n$ is BP consistent with $E_{0}$, then it is also Ziezold consistent with $E_{0}$. If $\hat{E}_{n}$ is Ziezold consistent with $E_{0}$, then $\hat{E}_{n}$ is also BP consistent to $E_{0}$ when Condition~\ref{C:Heine-Borel} holds and $\hat{E}_{n}$ is almost surely eventually bounded. 
\end{lemma}

\begin{proof}[of Lemma~\ref{lem:Z-SC-BP-SC-WC}]
To show the first part of Lemma~\ref{lem:Z-SC-BP-SC-WC}, it suffices to show that, for an arbitrary \(\omega \in \Omega\), the condition \(\lim_{n \rightarrow \infty} \sup_{m \in \hat{E}_n(\omega)} d(m, E_0) = 0\) implies
\[
\bigcap_{n=1}^{\infty} \overline{\bigcup_{k=n}^{\infty} \hat{E}_k(\omega)} \subset E_0.
\]
Suppose, for some \(\omega \in \Omega\), we have \(\lim_{n \rightarrow \infty} \sup_{m \in \hat{E}_n(\omega)} d(m, E_0) = 0\), but 
\[
\bigcap_{n=1}^{\infty} \overline{\bigcup_{k=n}^{\infty} \hat{E}_k(\omega)} \not\subset E_0.
\]
Then, there exists a subsequence \(e_{n_k} \in \hat{E}_{n_k}(\omega)\) and a point \(e \notin E_0\) such that \(d(e_{n_k}, e) \to 0\) as \(k \to \infty\). This implies
\[
\lim_{k \to \infty} \sup_{m \in \hat{E}_{n_k}(\omega)} d(m, E_0) \ge \lim_{k \to \infty} d(e_{n_k}, E_0) \ge \lim_{k \to \infty} \{ d(e, E_0) - d(e_{n_k}, e) \} = d(e, E_0) > 0,
\]
which contradicts the assumption that \(\sup_{m \in \hat{E}_n(\omega)} d(m, E_0) \to 0\). This shows the first part.

To show the second part, assume that Condition~\ref{C:Heine-Borel} holds. It suffices to show that if \(\hat{E}_n(\omega)\) is eventually bounded and \(\bigcap_{n=1}^{\infty} \overline{\bigcup_{k=n}^{\infty} \hat{E}_k(\omega)} \subset E_0\), then
\[
\lim_{n \to \infty} \sup_{m \in \hat{E}_n(\omega)} d(m, E_0) = 0.
\]
Suppose the contrary: there exists \(\epsilon > 0\) and a subsequence \(e_{n_k} \in \hat{E}_{n_k}(\omega)\) such that \(d(e_{n_k}, E_0) \ge \epsilon\) for all \(k\). Since \(\hat{E}_n(\omega)\) is eventually bounded and \(M\) satisfies the Heine–Borel property, the sequence \(\{e_{n_k}\}\) has a convergent subsequence \(e_{n_{k_m}} \to e\) for some \(e \in M\). Moreover, the assumption \(\bigcap_{n=1}^{\infty} \overline{\bigcup_{k=n}^{\infty} \hat{E}_k(\omega)} \subset E_0\) implies that \(e \in E_0\). Thus,
\[
d(e_{n_{k_m}}, E_0) \le d(e_{n_{k_m}}, e) \to 0
\] as $m\rightarrow \infty$,
which contradicts \(d(e_{n_k}, E_0) \ge \epsilon\) for all $k\ge 1$. Hence, we conclude that
\[
\lim_{n \to \infty} \sup_{m \in \hat{E}_n(\omega)} d(m, E_0) = 0,
\]
as desired.
\end{proof}

As discussed in the main text, Ziezold consistency holds under Conditions~\ref{C:Kuratowski}, \ref{C:separable}, \ref{C:rho-continuous}, and \ref{C:locally-integrable}. We formalize this observation in the following lemma.

\begin{lemma}\label{lem:Z-SC}
Assume Conditions~\ref{C:Kuratowski}, \ref{C:separable}, \ref{C:rho-continuous}, and \ref{C:locally-integrable}. Then, $\hat{E}_{n}$ is Ziezold consistent with $E_{0}$.
\end{lemma}

To show Lemma~\ref{lem:Z-SC}, we recall the notion of continuous convergence, which plays a central role in the argument. We also state two auxiliary lemmas that will be used in the proof.

\begin{definition}[Continuous convergence \citep{beer1985convergence}]\label{def: cont converge}
Let $f_n : M \rightarrow \mathbb{R}$ and $f : M \rightarrow \mathbb{R}$ be functions. We say that $f_n$ \emph{continuously converges} to $f$ if, for every $m \in M$ and every sequence $m_n \to m$, we have $\lim_{n \to \infty} f_n(m_n) = f(m)$.
\end{definition}

The following lemma shows that, under Condition~\ref{C:Kuratowski}, Ziezold consistency is implied by the almost sure continuous convergence of the empirical cost functions $F_n$ to the population cost function $F$.

\begin{lemma}\label{lem:epi-convergence-to-Z-SC}
Assume Condition~\ref{C:Kuratowski}. If $F_n$ continuously converges to $F$ almost surely, then $\hat{E}_n$ is Ziezold consistent with $E_0$.
\end{lemma}

\begin{proof}[of Lemma~\ref{lem:epi-convergence-to-Z-SC}]
Let \( l_0 = \inf_{m \in M_0} F(m) \) and \(\hat{l}_n(\omega) = \inf_{m \in M_n} F_n(\omega)(m)\). Fix \(\omega \in \Omega\) such that \(F_n(\omega)\) continuously converges to \(F\), and \(M_n(\omega)\) converges to \(M_0\) in the sense of Kuratowski.

We first show that \(\limsup_{n \to \infty} \hat{l}_n(\omega) \le l_0\). For any \(m_0 \in M_0\), the Kuratowski convergence of \(M_n\) to \(M_0\) guarantees the existence of a sequence \(m_n \in M_n(\omega)\) such that \(m_n \to m_0\). By the continuous convergence of \(F_n(\omega)\), we have \(F_n(\omega)(m_n) \to F(m_0)\). Since \(\hat{l}_n(\omega) \le F_n(\omega)(m_n)\), it follows that
\[
\limsup_{n \to \infty} \hat{l}_n(\omega) \le \limsup_{n \to \infty} F_n(\omega)(m_n) = F(m_0).
\]
As \(m_0 \in M_{0}\) was arbitrary, we conclude that \(\limsup_{n \to \infty} \hat{l}_n(\omega) \le l_0\).

To complete the proof, we show that \(\bigcap_{n=1}^{\infty} \overline{\bigcup_{k=n}^{\infty} \hat{E}_k(\omega)} \subset E_0\). The result is trivial if the intersection is empty, so assume \(\bigcap_{n=1}^{\infty} \overline{\bigcup_{k=n}^{\infty} \hat{E}_k(\omega)} \neq \emptyset\). Let \(b\) be an arbitrary point in this set. Then, there exists a subsequence \(m_{n_k}(\omega) \in \hat{E}_{n_k}(\omega)\) such that \(m_{n_k} \to b\). Kuratowski convergence implies \(b \in M_0\), and thus \(l_0 \le F(b)\).

From the definition of \(\hat{l}_n\) and $\hat{E}_{n}$, we have \(\hat{l}_{n_k}(\omega) \le F_{n_k}(\omega)(m_{n_k}) \le \hat{l}_{n_{k}}(\omega) +\varepsilon_{n_{k}}\). Taking limits yields
\[
\limsup_{n \to \infty} \hat{l}_n(\omega) \ge \limsup_{k \to \infty} \hat{l}_{n_k}(\omega) \ge \lim_{k \to \infty} F_{n_k}(\omega)(m_{n_k}) = F(b).
\]
Combining this with the earlier inequality \(\limsup_{n \to \infty} \hat{l}_n(\omega) \le l_0\), we conclude that \(F(b) = l_0\), and hence \(b \in E_0\). Since \(b\) was arbitrary, the claim follows.
\end{proof}

Next, we verify that Conditions~\ref{C:separable}, \ref{C:rho-continuous}, and \ref{C:locally-integrable} are sufficient to ensure the almost sure continuous convergence of $F_n$ to $F$. As mentioned in remark~\ref{rmk:transfer-ergodic}, guaranteeing the Ziezold-consistency \eqref{eq: Z-SC} and continuous convergence of $F_{n}$ to $F$ without specific condition or form of the cost function is the main technical difficulty in the proof of Theorem~\ref{thm: BP SC}. 

\begin{lemma}\label{lem:epi-convergence}
Assume Conditions~\ref{C:separable},\ref{C:rho-continuous}, and \ref{C:locally-integrable}. Then $F_n$ continuously converges to $F$ almost surely.
\end{lemma}

Before proceeding, we first establish some properties of $\pi_{m,r}(t)$ and $\Pi_{m,r}(t)$ in  Lemma~\ref{lem: scalarization properties} below.

\begin{lemma}\label{lem: scalarization properties} 
Recall that for a radius $r > 0$ and $m \in M$,  $\pi_{m,r}: T \rightarrow \mathbb{R}$ and $\Pi_{m,r}: T \rightarrow \mathbb{R}$ are defined as 
$\pi_{m,r}(t)= \inf_{m' \in B(m, r)} \cost(t, m')$, $\Pi_{m,r}(t)= \sup_{m' \in B(m, r)} \cost(t, m')$. 

(i) For all $m\in M$ and $t\in T$, $\pi_{m,r}(t)$ increases (and $\Pi_{m,r}(t)$ decreases) as $r$ decreases. (ii) Under Condition~\ref{C:rho-continuous}, $\pi_{m_{i},r_{i}}(t) \rightarrow \cost(t,m)$ and $\Pi_{m_{i},r_{i}}(t) \rightarrow \cost(t,m)$ as $m_{i}\rightarrow m$ and $r_{i}\rightarrow 0$ for all $t\in T$. (iii) Under Conditions~\ref{C:separable} and \ref{C:rho-continuous}, $\pi_{m,r}(t)$ and $\Pi_{m,r}(t)$ are measurable functions \citep[Theorem 3.4.(b)]{korf2001random}.
\end{lemma}

The first and second parts of Lemma~\ref{lem: scalarization properties} are immediate from the definition of $\pi_{m,r}$ and $\Pi_{m,r}$. For the third part, we refer to \cite{korf2001random}.

\begin{proof}[of Lemma~\ref{lem:epi-convergence}]
Since \(M\) is a separable metric space (Condition~\ref{C:separable}), it admits a countable dense subset \(D \subset M\). For each \(m \in D\), Condition~\ref{C:locally-integrable} ensures the existence of \(r_m > 0\) such that \(\pi_{m,r}(X)\) and \(\Pi_{m,r}(X)\) are integrable for all \(0 < r \le r_m\).

Fix \(m \in D\) and \(0 < r \le r_m\), and define the set
\[
\mathcal{A}_{m,r} = \left\{ \omega \in \Omega : \frac{1}{n} \sum_{i=1}^n \pi_{m,r}(X_i) \to E\{\pi_{m,r}(X)\} \text{ and } \frac{1}{n} \sum_{i=1}^n \Pi_{m,r}(X_i) \to E\{\Pi_{m,r}(X)\} \right\}.
\]
By the strong law of large numbers \citep[Theorem 2.4.1]{durrett2019probability}, each \(\mathcal{A}_{m,r}\) has probability one. Now define
\[
\mathcal{A} = \bigcap_{m \in D} \bigcap_{r \in \mathbb{Q} \cap (0, r_m]} \mathcal{A}_{m,r},
\]
which is also a probability-one subset of \(\Omega\), since it is a countable intersection of such sets.

Fix \(\omega \in \mathcal{A}\) that also satisfies Condition~\ref{C:Kuratowski}. We aim to show that \(F_n(\omega)\) continuously converges to \(F\), i.e., for any \(m_0 \in M_0\) and any sequence \(m_n \in M_n(\omega)\) such that \(m_n \to m_0\), we have \(F_n(\omega)(m_n) \to F(m_0)\).

Since \(D\) is dense in \(M\), one can choose a sequence \(m^l \in D\) such that \(m^l \to m_0\). For each \(l\), choose \(0 < r^l < r_{m^l}\) such that \(r^l \downarrow 0\), and
\[
m_0 \in \cdots \subset B(m^2, r^2) \subset B(m^1, r^1).
\]
Because \(m_n \to m_0\), for each \(l\) there exists \(N_l\) such that \(m_n \in B(m^l, r^l)\) for all \(n = N_l, N_{l}+1, \ldots \). Then,
\begin{align*}
    E\{\pi_{m^l, r^l}(X)\} &= \liminf_{n \to \infty} \frac{1}{n} \sum_{i=1}^n \pi_{m^l, r^l}(X_i) \le \liminf_{n \to \infty} F_n(\omega)(m_n), \\
    \limsup_{n \to \infty} F_n(\omega)(m_n) &\le \limsup_{n \to \infty} \frac{1}{n} \sum_{i=1}^n \Pi_{m^l, r^l}(X_i) = E\{\Pi_{m^l, r^l}(X)\}.
\end{align*}

From Condition~\ref{C:rho-continuous}, the function \(\cost(t, m)\) is continuous in \(m\) for each fixed \(t\). Therefore, \(\pi_{m^l, r^l}(t) \uparrow \cost(t, m_0)\) and \(\Pi_{m^l, r^l}(t) \downarrow \cost(t, m_0)\) as \(l \to \infty\) for all \(t \in T\). 

By the Monotone Convergence Theorem \citep[Theorem 1.6.6]{durrett2019probability}, we have
\[
E\{\pi_{m^l, r^l}(X)\} \to E\{\cost(X, m_0)\}.
\]
Similarly, since \(\Pi_{m^l, r^l}(X)\) is bounded above and below by integrable functions, the Dominated Convergence Theorem \citep[Theorem 1.6.7]{durrett2019probability} gives
\[
E\{\Pi_{m^l, r^l}(X)\} \to E\{\cost(X, m_0)\}.
\]
Combining the bounds, we obtain
\[
\lim_{n \to \infty} F_n(\omega)(m_n) = E\{\cost(X, m_0)\} = F(m_0).
\]
Since \(m_0\) and the sequence \(m_n \to m_0\) were arbitrary, this proves that \(F_n(\omega)\) continuously converges to \(F\).
\end{proof}

Conclusion of Lemma~\ref{lem:Z-SC} immediately follows from Lemmas~\ref{lem:epi-convergence-to-Z-SC} and \ref{lem:epi-convergence}.
\begin{proof}[of Lemma~\ref{lem:Z-SC}]
By Lemma~\ref{lem:epi-convergence}, $F_{n}$ continuously converges to $F$ almost surely under Conditions \ref{C:separable}, \ref{C:rho-continuous}, and \ref{C:locally-integrable}. Lemma~\ref{lem:epi-convergence-to-Z-SC} implies Ziezold consistency of $\hat{E}_{n}$ to $E_{0}$ with Condition~\ref{C:Kuratowski}.
\end{proof}

Combining the technical results Lemmas~\ref{lem:Z-SC-BP-SC-WC}, \ref{lem:Z-SC}, we immediately obtain Theorem~\ref{thm: BP SC}. 
\begin{proof}[of Theorem~\ref{thm: BP SC}]
By Lemma \ref{lem:Z-SC}, $\hat{E}_{n}$ is Ziezold consistent with $E_{0}$ under Conditions~\ref{C:Kuratowski}, \ref{C:separable}, \ref{C:rho-continuous}, and \ref{C:locally-integrable}. If Condition~\ref{C:Heine-Borel} and the almost sure eventual boundedness of $\hat{E}_{n}$ is additionally satisfied, then $\hat{E}_{n}$ is BP consistent with $E_{0}$ by Lemma~\ref{lem:Z-SC-BP-SC-WC}.
\end{proof}

\subsection{Technical details for Section~\ref{seC:various examples}}\label{app: technical details for Sec 2.4}

In this subsection, we present detailed proofs of Corollaries~\ref{cor: BP SC huckemann generalization} and \ref{cor: BP SC Hd}, relying on the result established in Theorem~\ref{thm: BP SC}. We also introduce several auxiliary results essential for completing these proofs.

To show Corollary~\ref{cor: BP SC huckemann generalization}, we establish the following two lemmas, which ensure Conditions~\ref{C:rho-continuous}, \ref{C:locally-integrable} and the almost sure eventual boundedness of $\hat{E}_{n}$ under assumptions in Corollary~\ref{cor: BP SC huckemann generalization}.

\begin{lemma}\label{lem: locally integrable}
Assume that $\cost=G\circ \rho$ for equicontinuous $\rho: T\times M \rightarrow [0,\infty)$ and continuous, non-decreasing, and $b$-subadditive $G: [0,\infty) \rightarrow [0,\infty)$. Then Conditions~\ref{C:rho-continuous} and \ref{C:locally-integrable} holds.
\end{lemma} 

\begin{proof}[of Lemma~\ref{lem: locally integrable}]
We first verify Condition~\ref{C:rho-continuous}. Since $G$ is continuous and $\rho(t,\cdot): M\rightarrow [0,\infty)$ is continuous for $t\in T$, $\cost(t,\cdot)=G\circ \rho(t,\cdot)$ is continuous for all $t\in T$. 

Secondly, we show that Condition~\ref{C:locally-integrable} holds. We point out that $G(x+y)\le bG(x)+bG(y)$ (for some $b\ge 1$) when $G$ is $b$-subadditive and non-decreasing; see Lemma~\ref{lem: H properties}. Since $\cost=\rho^{2} \ge 0$ and $\cost(X,m)$ is integrable, $0\le \pi_{m,r}(X) \le \cost(X,m)$ is integrable for all $r>0$. 

To show the integrability of $\Pi_{m,r}(X)$, we choose $\delta>0$ so that $\lvert  \rho(t,m')- \rho(t,m) \rvert <1$ for all $m'\in B(m,\delta)$ and for all $t\in T$. Such $\delta$ exists by equicontinuity of $\rho$. Then 
\begin{equation*}
\begin{split}
 \cost(X,m')& =G\{ \rho(X,m')\} \le G\{  \rho(X,m)+\mid \rho(X,m)-\rho(X,m') \mid \}  \\
& \le b G \{\rho(X,m)\} + bG \{\mid \rho(X,m)-\rho(X,m') \mid )\}\\ 
&\le b \mathfrak{c}(X,m) + b G(1)    
\end{split}    
\end{equation*} 
 for $m'\in B(m,\delta)$. Thus $\Pi_{m,r}(X)$ is integrable for all $0<r<\delta$, and Condition~\ref{C:locally-integrable} holds.
\end{proof}

\begin{lemma}\label{lem: coercive} Let $\cost=G\circ \rho$ for equicontinuous, coercive $\rho: T\times M \rightarrow [0,\infty)$ and 
continuous, non decreasing, $b$-subadditive (for some $b\ge 1$) $G: [0,\infty) \rightarrow [0,\infty)$ that satisfies $\lim_{x\rightarrow \infty} G(x) = \infty$. Assume Conditions~\ref{C:Kuratowski} and \ref{C:separable}. Then, $\hat{E}_{n}$ is eventually bounded almost surely.
\end{lemma}

\begin{proof}[of Lemma~\ref{lem: coercive}]
By Lemma~\ref{lem: locally integrable}, Conditions~\ref{C:rho-continuous} and \ref{C:locally-integrable} holds. Thus, by Theorem~\ref{thm: convergence of minimum}, 
\begin{equation}
    \label{eq:S7-1}
\limsup_{n\rightarrow \infty} \hat{l}_{n}(\omega) \le l_{0}
\end{equation}
almost surely, where $\hat{l}_{n}(\omega) = \inf_{m\in M_{n}} F_{n}(\omega)(m)$ and ${l}_{0} = \inf_{m\in M_{0}} F(m)$.

Since $\rho$ is coercive, there exist $o \in M$ and $C>0$ that satisfy the two requirements of the coercive function, as stated in Definition~\ref{def: equicontinuous} (ii). Let $k(n)(\omega)= | \{i : i = 1,\ldots, n \ {\rm and } \  \rho(X_{i}(\omega),o)<C\} |$ for $n=1,2,\ldots $ and $\omega \in \Omega$.
By the strong law of large numbers, we have  
\begin{equation}
    \label{eq:S7-2}
    \frac{k(n)(\omega)}{n}\rightarrow  P[\{\rho(X,o)<C\}]>0
\end{equation}
for almost all $\omega\in \Omega$.
Fix an $\omega \in \Omega$ that satisfies both (\ref{eq:S7-1}) and (\ref{eq:S7-2}).

Note that the set of such $\omega$ has probability $1$. Now it is enough to show that for such $\omega\in \Omega$, $$r_{n}(\omega)=\sup \{d( m_{n},o) :  m_{n}\in \epsilon_{n}\text{-}\arg\min_{m\in M_{n}} F_{n}(\omega)(m)\}$$ is bounded. Suppose $r_{n}(\omega)$ is not bounded above, so that there exists a subsequence $r_{n_{a}}(\omega)\rightarrow \infty$ as $a\rightarrow \infty$. Then we can choose a sequence $m_{n_{a}}(\omega)\in  \epsilon_{n_{a}}\text{-}\arg\min_{m\in M_{n_{a}}} F_{n_{a}}(\omega)(m)$ such that $d\{m_{n_{a}}(\omega), o\} \to \infty$. Let $A_{n}$ be a real sequence satisfying $A_{n}\rightarrow \infty$ such that $\rho(t,m_{n})>A_{n}$ for all $t\in T$ satisfying $\rho(t,o)<C$.
Then, we obtain 
\begin{equation*}
\begin{split}
\epsilon_{n_{a}}+\hat{l}_{n_{a}}(\omega)
&  \ge F_{n_{a}}(\omega)\{m_{n_{a}}(\omega)\}=\frac{1}{n_{a}}\sum_{i=1}^{n_{a}} G \circ \rho \{X_{i}(\omega),m_{n_{a}}(\omega)\}  \\
& \ge \frac{1}{n_{a}}\sum_{i \in K(n_a)(\omega)} G\circ\rho \{ X_{i}(\omega),m_{n_{a}}(\omega) \} \\
&\ge \frac{k(n_{a})(\omega)}{n_{a}} G\{A_{n_{a}}(\omega)\}.
\end{split}
\end{equation*} 
Collecting (\ref{eq:S7-2}) and the facts that  $G(x)\rightarrow \infty$ as $x\rightarrow\infty$, $\epsilon_{n}\rightarrow 0$, and $A_{n}\rightarrow \infty$, we have  $\hat{l}_{n_{a}}(\omega) \rightarrow \infty $ as $a\rightarrow \infty$. It is a contradiction to ${\limsup }_{n\rightarrow \infty} \hat{l}_{n}(\omega)\le l_{0}$. 

Thus, $r_{n}(\omega)$ is bounded above almost surely, $\epsilon_{n}\text{-}\arg\min_{m\in M_{n}} F_{n}(\omega)(m)$ is eventually bounded almost surely.
\end{proof}

Corollary~\ref{cor: BP SC huckemann generalization} is immediate by combining the results of Theorem~\ref{thm: BP SC} and Lemmas~\ref{lem: locally integrable} and \ref{lem: coercive}.

\begin{proof}[of Corollary~\ref{cor: BP SC huckemann generalization}]
By Lemmas~\ref{lem: locally integrable} and \ref{lem: coercive}, Conditions~\ref{C:rho-continuous} and \ref{C:locally-integrable} hold, and $\hat{E}_{n}$ is eventually bounded almost surely. Thus, Theorem~\ref{thm: BP SC} gives that $\hat{E}_{n}$ is BP consistent with $E_{0}$.
\end{proof}

As referenced in the main text, Corollary~\ref{cor: BP SC Hd} is an immediate consequence of Corollary~\ref{cor: BP SC huckemann generalization}.

\begin{proof}
It is enough to show $d: T\times M\rightarrow \mathbb{R}$ satisfies equicontinuity and coercivity, when $T=M$. 

We first show equicontinuity of $d$. For all $t\in T$ and $m,m'\in M$, $\mid d(t,m')- d(t,m) \mid \le d(m',m)$ holds by tringle inequality. Thus, Equicontinuity is trivially satisfied by setting $\delta({\epsilon, m}) = \epsilon$.

The coercivity of $d$ can be shown as below. Choose an arbitrary $o\in O$, and $C>0$ such that $\text{pr}\{d(X,o)<C\}>0$. Such $C$ exists since $\lim_{n\rightarrow \infty} \text{pr}\{d(X,o)<n\}=1$ for all $o\in M$. Then, for arbitrary sequence $m_{n}\in M$ with $d(o,m_{n})\rightarrow \infty$, let $A_{n} = d(o,m_{n})-C$. Then, $A_{n}\rightarrow \infty$ and 
$$d(t,m_{n})\ge d(o,m_{n})-d(t,o)> A_{n}$$ for all $t\in T$ satisfying $d(t,o)<C$, by triangle inequality of $d$.
\end{proof}



\subsection{Technical details for Section~\ref{sec:all_applications}}\label{app: technical details for Sec 2.5}

This section provides proof for Theorem~\ref{thm: BP consistency of PGA} and several intermediate results on the principal geodesic analysis in hyperspheres. 
Recall that the population and empirical $k$-dimensional principal geodesic submanifolds (or subspheres), denoted by $S^{k}$ and $\hat{S}^{k}$ respectively, are the population and empirical generalized \fre means. For such an identification, we denote the $(k-1)$-dimensional principal geodesic subspheres, obtained in the previous stage of {forward dimension reduction}, by $S^{k-1}$ and $\hat{S}^{k-1}$. 

The data space $T$ is the metric space $(\mathbb{S}^{p-1}, d)$, while the parameter space $M$ is $(\mathfrak{S}^{k},h_{d})$,
where $\mathfrak{S}^{k}$ is the set of all $k$-dimensional subspheres in $\mathbb{S}^{p-1}$, endowed with geodesic distance $d$ and $h_{d}$ is the Hausdorff distance induced by $d$. 
We set 
\begin{equation*}
     M_0 =\mathfrak{S}^{k}({S}^{k-1}),  \ 
     M_n = \mathfrak{S}^{k}(\hat{S}^{k-1}),  \ \text{and} \ \cost(t,m)=d^{2}(t,m)
\end{equation*} 
where, for a set $S \subset \mathbb{S}^{p-1}$, $\mathfrak{S}^{k}({S})$ is the set of all $k$-dimensional subspheres in $\mathbb{S}^{p-1}$ containing $S$.


Before demonstrating the consistency of $\hat{S}^{k}$, we first establish the existence of the population and empirical principal geodesic subspheres, $S^{k}$ and $\hat{S}^{k}$. While this result, stated in Theorem~\ref{thm: existence of principal great sphere} below, is of independent interest, the intermediate results (Lemmas~\ref{lem: continuous rho} and \ref{lem: compact domain}) obtained while verifying Theorem~\ref{thm: existence of principal great sphere} will be utilized in the proof of Theorem~\ref{thm: BP consistency of PGA}.

In Lemma~\ref{lem: continuous rho} below, the map $d: \mathbb{S}^{p-1} \times \mathfrak{S}^{k}\rightarrow [0,\infty) $ is given by $d(x, S) = \inf\{d(x,y) : y \in S\}$.


\begin{lemma}\label{lem: continuous rho}
The map 
$d: \mathbb{S}^{p-1} \times \mathfrak{S}^{k}\rightarrow [0,\infty) $ is equicontinuous in the sense of Definition~\ref{def: equicontinuous}. (i), with respect to the Hausdorff distance $h_{d}$. Specifically, we have $$|d(s,A)-d(s,B)|\le h_{d}(A,B)$$ for any $A, B \in \mathfrak{S}^{k}$ and $s\in \mathbb{S}^{p-1}$.
Moreover, for any distribution $P$ on $\mathbb{S}^{p-1}$, $S\rightarrow E\{d^{2}(X,S)\}$ is a continuous function of $S$ on the domain $(\mathfrak{S}^{k}, h_{d})$.
\end{lemma}

\begin{proof}[of Lemma~\ref{lem: continuous rho}]
By symmetry, it is enough to show $d(s, A)\le d(s, B) + h_{d}(A ,B)$ only. For arbitrary $a\in A$ and $b\in B$, 
\begin{equation}
   d(s,a)\le d(s,b)+ d(a,b) \label{eq:aA}
\end{equation} holds. Taking the infimum over $a \in A$ for both sides of (\ref{eq:aA}), we get $d(s,A)\le d(s,b) + \inf_{a\in A} d(a,b)$. Since $\inf_{a\in A} d(a,b) \le \sup_{b\in B}\inf_{a\in A} d(a,b)$, $\inf_{a\in A} d(a,b) \le h_{d}(A,B)$.
Thus, we obtain 
\begin{equation}
   d(s,A)\le d(s,b)+ h_{d}(A,B). \label{eq:bB}
\end{equation}
Taking the infimum over $b\in B$ for both sides of (\ref{eq:bB}), we get $d(s,A)\le d(s,B)+ h_{d}(A,B)$, as required.

The second part of Lemma~\ref{lem: continuous rho} can be shown by using the first part of Lemma~\ref{lem: continuous rho}. For any $S_{1}, S_{2} \in \mathfrak{S}^{k}$, 
\begin{equation*}
\begin{split}
|E \{d^{2}(X,S_{1})\}- E \{ d^{2}(X,S_{1}) \} | & \le E \{|d^{2}(X,S_{1})-d^{2}(X,S_{2}) |\} \\
& = E \{ |d(X,S_{1})-d(X,S_{2}) |(d(X,S_{1})+d(X,S_{2})) \} \\
& \le 2\pi E\{ |d(X,S_{1})-d(X,S_{2}) | \}\\
&\le 2\pi h_{d}(S_{1},S_{2}) 
\end{split}    
\end{equation*}
holds. (Here, $\pi \approx 3.14159$ is the mathematical constant.) The first part of Lemma~\ref{lem: continuous rho} was used in the last inequality. Note that geodesic distance $d$ on $\mathbb{S}^{p}$ is bounded above by $\pi$.
\end{proof}

Since $\mathbb{S}^{p-1}$ is the unit sphere embedded in $\mathbb{R}^p$, any $k$-dimensional geodesic submanifold is the intersection of $\mathbb{S}^{p-1}$ and a $k$-dimensional subspace of $\mathbb{R}^p$. Thus, a set $V_k = \{v_{1}, v_{2}, \ldots, v_{k}\}$ of orthonormal vectors in $\mathbb{R}^{p}$ defines a $k$-dimensional geodesic submanifold (or subsphere), which we denote by $S(V_k)$ or $S(v_{1}, \ldots, v_{k})$. That is,  
$$S(V_k) = S(v_{1}, \ldots, v_{k}) = \mathbb{S}^{p-1}\cap \text{span}\{v_{1}, v_{2}, \ldots, v_{k}\}.$$
We also write $\text{span}(V_k)$ for $\text{span}\{v_{1}, v_{2}, \ldots, v_{k}\}$. 



\begin{lemma}\label{lem: compact domain} For $k=1,\ldots, p-1$, the following hold.

(i) $\mathfrak{S}^{k}$ is compact with respect to $h_{d}$. 

(ii) For any orthonormal subset $V_{k}=\{v_{1}, v_{2}, \ldots, v_{k}\}$ of $\mathbb{R}^{p}$, $\mathfrak{S}^{k}\{S(V_{k})\}$ is a closed subset of $\mathfrak{S}^{k}$ and is therefore compact with respect to $h_{d}$. 

(iii) For orthonormal subsets $V_k = \{v_{1}, v_{2}, \ldots, v_{k}\}$ and $V_{n,k} = \{v_{n,1}, \ldots, v_{n,k} \}$ of $\mathbb{R}^{p}$, suppose that  $d(v_{n,i}, v_{i})\rightarrow 0$ as $n\rightarrow \infty$ for $i=1,\ldots, k$. Then, $\mathfrak{S}^{k}\{S(V_{n,k})\}$ converges to $\mathfrak{S}^{k}\{S(V_k)\}$ in the sense of Kuratowski.
\end{lemma}

To show Lemma~\ref{lem: compact domain}, we bound the Hausdorff distance between two geodesic submanifolds $S$ and $S'$ by the distance between two direction vectors, each of which lies in $S$ and $S'$, respectively.

\begin{lemma}\label{lem: hausdorff distance between great sphere}
For a $k=1,\ldots, p-1$, let $V_k = \{v_{1}, \ldots, v_{k}\}$ and $V_k' = \{v'_{1}, \ldots, v'_{k}\}$ be two orthonormal subsets of $\mathbb{R}^{p}$, satisfying $0\le d(v_{i}, v_{i}')\le \cos^{-1}(\frac{k-1}{k})$ for all $i=1,2,\ldots, k$. Then for $S=S(v_{1}, \ldots, v_{k})$ and $S'=S(v'_{1}, \ldots, v'_{k})$, $$h_{d}(S,S')\le \cos^{-1} [1-k  + k \cos\{ \max_{i=1,2,\ldots, k}d(v_{i}, v_{i}') \}].$$
Thus, for orthonormal subsets $V_k = \{v_{1}, v_{2}, \ldots, v_{k}\}$ and $V_{n,k} = \{v_{n,1}, \ldots, v_{n,k} \}$ of $\mathbb{R}^{p}$, if $d(v_{i}, v_{n,i})\rightarrow 0$ as $n\rightarrow \infty$ for $i=1,2,\ldots, k$, $h_{d}\{S(V_k), S(V_{n,k})\}\rightarrow 0$ as $n\rightarrow \infty$.
\end{lemma}

\begin{proof}[of Lemma~\ref{lem: hausdorff distance between great sphere}]
For any $s\in S$, there exists $x=(x_{1}, \ldots, x_{k})\in \mathbb{S}^{k-1}$ such that 
$s=x_{1}v_{1}+\cdots + x_{k}v_{k}$. For such $x$, $x_{1}v'_{1}+\cdots + x_{k}v'_{k}\in S'$. Then,  

\begin{align}
&\inf_{s'\in S'} d(s,s') \le d(x_{1}v_{1}+\cdots + x_{k}v_{k}, x_{1}v'_{1}+\cdots + x_{k}v'_{k}) \nonumber  \\
& = \cos^{-1} \{(x_{1}v_{1}+\cdots + x_{k}v_{k})^{T} (x_{1}v_{1}'+\cdots + x_{k}v_{k}') \} \label{eq: angle and norm 1} \\ 
& = \cos^{-1} \{1-\|x_{1}(v_{1}-v_{1}')+\cdots + x_{k}(v_{k} -v_{k}')\|^{2} /2\}.\label{eq: angle and norm 2}
\end{align}

\eqref{eq: angle and norm 1} follows from the definition of geodesic distance $d$ in $\mathbb{S}^{p-1}$, and \eqref{eq: angle and norm 2} follows from 
\begin{equation*}
\begin{split}
&\|x_{1}(v_{1}-v_{1}')+\cdots + x_{k}(v_{k} -v_{k}')\|^{2}  \\
& = \|x_{1}v_{1}+\cdots + x_{k}v_{k}\|^{2} + \| x_{1}v_{1}'+\cdots + x_{k}v_{k}'\|^{2} - 2 (x_{1}v_{1}+\cdots + x_{k}v_{k})^{T} (x_{1}v_{1}'+\cdots + x_{k}v_{k}') \\ 
& = 2-2(x_{1}v_{1}+\cdots + x_{k}v_{k})^{T} (x_{1}v_{1}'+\cdots + x_{k}v_{k}').
\end{split}    
\end{equation*}

Using the property of Euclidean norm and the Cauchy-Schwarz inequality, we get 
\begin{align}
& \|x_{1}(v_{1}-v_{1}')+\cdots + x_{k}(v_{k} -v_{k}')\| \le |x_{1}| \|v_{1}-v_{1}'\| + \cdots + |x_{k}| \|v_{k}-v_{k}'\| ~ \nonumber  \\
& \le \sqrt{x_{1}^{2} + \ldots + x_{k}^{2}} \sqrt{ \|v_{1}-v_{1}'\|^{2} + \cdots +  \|v_{k}-v_{k}'\|^{2}} ~ \nonumber \\
& \le \sqrt{k} \max_{i=1,\ldots, k}{ \|v_{i}-v_{i}'\|} \label{eq: angle and norm 3}.
\end{align}

Combining \eqref{eq: angle and norm 2} and \eqref{eq: angle and norm 3}, and the fact that $\cos^{-1}$ is a decreasing function on the interval $[0,1]$, 
we get
\begin{align}
& \inf_{s'\in S'} d(s,s') \nonumber  \\
& \le \cos^{-1} (1-k  \max_{i=1,\ldots, k}{ \|v_{i}-v_{i}'\|^{2}} /2) \nonumber \\ 
& = \cos^{-1} [1-k  + k \cos\{ \max_{i=1,2,\ldots, k}d(v_{i}, v_{i}') \}] .\nonumber
\end{align}

Since $s\in S$ was arbitrary, we have $$\sup_{s\in S} \inf_{s'\in S'} d(s,s') \le \cos^{-1} [1-k  + k \cos\{ \max_{i=1,2,\ldots, k}d(v_{i}, v_{i}') \}].$$

By the symmetry, we also have $$\sup_{s'\in S'} \inf_{s\in S} d(s,s') \le \cos^{-1} [1-k  + k \cos\{ \max_{i=1,2,\ldots, k}d(v_{i}, v_{i}') \}].$$ The definition of Hausdorff distance follows the conclusion.
\end{proof}

\begin{proof}[of Lemma~\ref{lem: compact domain}]
Proof of (i). Let $S_{n}$ be an arbitrary $\mathfrak{S}^{k}$-valued sequence. It is enough to find a subsequence of $S_{n}$, whose limit $S_{0}$ is in $\mathfrak{S}^{k}$. Each $S_{n}$ can be represented by $S(V_{n,k+1})$ for an orthonormal subset $V_{n,k+1} = \{v_{n,1}, v_{n,2}, \ldots, v_{n,k}, v_{n,k+1}\}$ in $\mathbb{R}^{p}$. Since $\mathbb{S}^{p-1}$ is compact (with respect to $d$), there exists a subsequence $n_{l}$ and $v_{1}, \ldots, v_{k}, v_{k+1}\in {\mathbb{S}}^{p-1}$ such that $d(v_{n_{l}, i}, v_{i}) \rightarrow 0$ as $l\rightarrow \infty$ for $i=1,\ldots,k+1$. Notice that $V_{0}=\{v_{1}, v_{2}, \ldots, v_{k},v_{k+1} \}$ is an orthonormal subset of $\mathbb{R}^{p}$, and $S(V_{0}) \in \mathfrak{S}^{k}$. By Lemma~\ref{lem: hausdorff distance between great sphere}, $h_{d}\{S_{n_{l}}, S(V_{0})\}\rightarrow 0$ as $l\rightarrow \infty$. Thus we can conclude that $S_{0}=S(V_{0})\in \mathfrak{S}^{k}$. 

(ii) The second part can be shown similarly. 
Let $S_{n}$ be an arbitrary sequence in $\mathfrak{S}^{k}\{S(V_{k})\}$  satisfying $\lim_{n\rightarrow \infty} h_{d}(S_{n}, S_{0})=0$ for some $S_{0}\in \mathfrak{S}^{k}$. It is enough to show $S_{0}\in \mathfrak{S}^{k}\{S(V_{k})\}$.

To apply Lemma~\ref{lem: hausdorff distance between great sphere} to this situation, we represent, for each $n$, the geodesic subsphere $S_{n}$ by $S_{n}=S(v_{1}, v_{2}, \ldots, v_{k}, v_{n,k+1})$ for an $v_{n,k+1}\in \mathbb{S}^{p-1}\cap \{\text{span}(V_{k})\}^{\perp}$. (Recall that $v_i$'s ($i=1,\ldots,k$) above are exactly the orthonormal vectors in $V_k  = \{v_{1}, v_{2}, \ldots, v_{k}\}$.) Since $\mathbb{S}^{p-1}\cap \{\text{span}(V_{k})\}^{\perp}$ is compact, there exist a subsequence $n_{m}$ and $v_{k+1}\in \mathbb{S}^{p-1}\cap \{\text{span}(V_{k})\}^{\perp}$ such that $\lim_{m\rightarrow \infty} d(v_{n_{m}, k+1}, v_{k+1})=0$. Thus, $h_{d}\{S_{n_{m}}, S(v_{1}, v_{2}, \ldots, v_{k}, v_{k+1}) \} \rightarrow 0$ as $m\rightarrow \infty$, by Lemma~\ref{lem: hausdorff distance between great sphere}.  Since $h_{d}(S_{n}, S_{0})\rightarrow 0$, we get $S_{0}=S(v_{1}, v_{2}, \ldots, v_{k}, v_{k+1}) \in  \mathfrak{S}^{k}\{S(V_{k})\}$.

(iii) 
To check the first condition of Kuratowski convergence of $\mathfrak{S}^{k}\{S(V_{n,k})\}$ to $\mathfrak{S}^{k}\{S(V_{k})\}$, 
let 
$S_{n} \in \mathfrak{S}^{k}\{S(V_{n,k})\}$ be an arbitrary sequence, and write, for each $n$, 
$S_{n}=S(v_{n,1}, v_{n,2}, \ldots, v_{n,k}, v_{n,k+1})$ for some $v_{n,k+1} \in \mathbb{S}^{p-1} \cap \{\text{span}(V_{n,k})\}^{\perp}$.
%
Since $\mathbb{S}^{p-1}$ is compact, there exist $v_{k+1}\in \mathbb{S}^{p-1}$ and subsequence $n_{m}$ such that $d(v_{n_{m},k+1}, v_{k+1})\rightarrow 0$ as $m\rightarrow \infty$. 
One can easily show that $v_{k+1}$ is orthogonal to $v_{1}, \ldots, v_{k}$, we have $S_{0} = S(v_{1}, v_{2}, \ldots, v_{k},v_{k+1}) \in \mathfrak{S}^{k}\{S(V_{k})\}$ and $h_{d} \{S_{n_{m}}, S_0\} \rightarrow 0$ by Lemma~\ref{lem: hausdorff distance between great sphere}.  


To check the {second} condition of Kuratowski convergence, 
choose an arbitrary $S_{0}^{k}\in \mathfrak{S}^{k}\{S(V_{k})\}$, which can be represented by $S_{0}^{k}=S(v_{1}, v_{2}, \ldots, v_{k}, v_{k+1})$ for some $v_{k+1}\in {\mathbb{S}}^{p-1}\cap  \{\text{span}(V_{k})\}^{\perp}$. 
For $$v_{n,k+1}=\frac{ v_{k+1}-\text{proj}_{\text{span}(V_{n,k})} v_{k+1}  }{ \|v_{k+1}-\text{proj}_{\text{span}(V_{n,k})} v_{k+1} \|_{2} },$$
let
$S_{n}^{k}= S(v_{n,1}, v_{n,2}, \ldots, v_{n, k}, v_{n,k+1})\in \mathfrak{S}^{k}(S_{n}^{k-1})$. Then, $h_{d} (S_{n}^{k}, S_{0}^{k})\rightarrow 0$ as $n \rightarrow \infty$, by Lemma~\ref{lem: hausdorff distance between great sphere}.

\end{proof}

Combining Lemmas~\ref{lem: continuous rho} and \ref{lem: compact domain}, we conclude that $S^{k}$ exists for any distribution $P$, since any continuous function has a minimizer on a compact set. The empirical estimators $\hat{{S}}^{k}$ can be regarded as a special case of $S^{k}$ derived from a distribution that assigns probability $1 / n$ to each fixed $X_1, X_2, \ldots, X_n$. Thus, for any distribution $P$ and for any observed $X_1, X_2, \ldots, X_n$, both $S^{k}$ and $\hat{S}^{k}$ exist for $k=0,1,\ldots, p-1$. This result is summarized below. 

\begin{theorem}\label{thm: existence of principal great sphere}
For any distribution $P$ on $\mathbb{S}^{p-1}$ and for $k=0,\ldots, p-1$, both ${S}^{k}$ and ${\hat{S}}^{k}$ exist.
\end{theorem} 

Finally, we state a converse of Lemma~\ref{lem: hausdorff distance between great sphere}, which will be used in our proof of Theorem~\ref{thm: BP consistency of PGA}.  


\begin{lemma}\label{lem: hausdorff distance between great sphere converse}
Fix a $k$ among $\{1,\ldots, p\}$, and let 
$$S_{0}^{0} \subset \cdots \subset S_{0}^{k-1},$$
where $S_{0}^{i}\in \mathfrak{S}^{i}$ ($i=0,\ldots,k-1$), be a nested sequence of geodesic subspheres. Likewise, for each $n$, let 
$$S_{n}^{0}\subset  \cdots \subset S_{n}^{k-1},$$
where $S_{n}^{i}\in \mathfrak{S}^{i}$, be artirary while satisfying $h_{d}(S_{0}^{i}, S_{n}^{i}) \rightarrow 0$ as $n\rightarrow \infty$.
%
%
Let $\{v_{0}, \ldots, v_{k-1}\}$ be an orthonormal subset of $\mathbb{R}^{p}$ satisfying that 
$$S_{0}^{0} = \{v_{0}\}, 
S_{0}^{1}=\mathbb{S}^{p-1}\cap \text{span}\{v_{0}, v_{1}\}, \ldots,
S_{0}^{k-1}=\mathbb{S}^{p-1}\cap \text{span}\{v_{0},  \ldots, v_{k-1}\}.$$ 
Similarly, for each $n$, let $\{v_{n,0},\ldots, v_{n,k}\}$ be an orthonormal subset of $\mathbb{R}^{p}$ satisfying that 
$S_{n}^{0} = \{v_{n,0}\}$ and $S_{n}^{i}=\mathbb{S}^{p-1}\cap \text{span}\{v_{n,0},  ..., v_{n,i}\}$ ($i=1,\ldots, k-1$). Then, for $i = 0,\ldots, k-1$,
$$\min\{d(v_{i}, v_{n,i}),d(v_{i}, -v_{n,i})\}\rightarrow 0 \quad \mbox{as} \ n \to \infty.$$ 
Adjusting the signs of $v_{n,i}$'s, we can say that $d(v_{i}, v_{n,i})\rightarrow 0$ as $n\rightarrow \infty$.
\end{lemma}

\begin{proof}[of Lemma~\ref{lem: hausdorff distance between great sphere converse}]
Proof of Lemma~\ref{lem: hausdorff distance between great sphere converse} is done by induction for $i$.

For $i=0$, since $h_{d}(S_{0}^{0}, S_{n}^{0}) = d(v_{0}, v_{n,0}) \rightarrow 0$ as $n \to \infty$, $d(v_{0}, v_{n,0}) \rightarrow 0$   as well. 

For $i\ge 1$, suppose $d(v_{j}, v_{n,j}) \rightarrow 0$ for $j=0,\ldots, i-1$ but $v_{n,i}$ does not converges to $v_{i}$. Since $\mathbb{S}^{p-1}$ is compact, there exist a subsequence of $v_{n_{m},i}$ and $v_{i}'\neq v_{i}\in \mathbb{S}^{p-1}$, such that $\lim_{m\rightarrow \infty} d(v_{n_{m},i}, v_{i}')=0$. Then, clearly $v_{i}' \perp v_{0}, \ldots, v_{i-1}$. Applying Lemma~\ref{lem: hausdorff distance between great sphere}, we get $$\lim_{m\rightarrow \infty} h_{d}(S_{n_{m}}^{i} , \mathbb{S}^{p-1}\cap \text{span} \{v_{0}, v_{1}, \ldots, v_{i-1}, v_{i}'\})=0.$$ Since $\lim_{n\rightarrow \infty}h_{d}(S_{n}^{i}, S_{0}^{i}) = 0$, we get $S_{0}^{i} = \mathbb{S}^{p-1}\cap \text{span} \{v_{0}, v_{1}, \ldots, v_{i-1}, v_{i}'\})$, thus $v_{i}'= - v_{i}$. This concludes that $\min\{d(v_{i}, v_{n,i}),d(-v_{i}, v_{n,i})\}\rightarrow 0$. 

By mathematical induction, we obtain $\min\{d(v_{i}, v_{n,i}),d(v_{i}, -v_{n,i})\}\rightarrow 0$ as $n\rightarrow \infty$ for all $i=0,\ldots, k$.
\end{proof}

Our general result in Theorem~\ref{thm: BP SC} and auxiliary results in Lemmas~\ref{lem: continuous rho}, \ref{lem: compact domain}, \ref{lem: hausdorff distance between great sphere} and \ref{lem: hausdorff distance between great sphere converse} are now used to show Theorem~\ref{thm: BP consistency of PGA} on the consistency of the empirical principal geodesic subspheres.  


\begin{proof}[of Theorem~\ref{thm: BP consistency of PGA}]
We will show the conclusion of Theorem~\ref{thm: BP consistency of PGA} by induction for $k$, applying Theorem~\ref{thm: BP SC}. 

Let $\{v_{0}, v_{1}, \ldots, v_{p-1}\}$ be an orthonormal basis of $\mathbb{R}^{p}$, satisfying $S^{0}= \{v_{0}\}$ and $S^{i}= \{v_{0}, v_{1}, \ldots, v_{i}\}$ for $i=1,\ldots, p-1$. We also write $\{\hat{v}_{0}, \hat{v}_{1}, \ldots, \hat{v}_{p-1}\}$ for an orthonormal basis of $\mathbb{R}^{p}$, satisfying $\hat{S}^{0}= \{\hat{v}_{0}\}$ and $\hat{S}^{i}= \{\hat{v}_{0}, \hat{v}_{1}, \ldots, \hat{v}_{i}\}$ for $i=1,\ldots, p-1$. (Note that $\{\hat{v}_{0}, \hat{v}_{1}, \ldots, \hat{v}_{p-1}\}$ is based on a sample of size $n$.)

For $k=0$, $S^{0} = \{v_0\}$ and $\hat{S}^{0} = \{\hat{v}_0\}$ are the population and sample \fre means, respectively. Therefore, the generalized \fre mean framework applies with 
$$ 
T=\mathbb{S}^{p-1}, (M,d)=(\mathbb{S}^{p-1},d), M_{0}=M_{n}=M, \ \mbox{and} \ \cost(t,m)=d^{2}(t,m)
$$
for $t\in T, m\in M$. To verify that $d(\hat{v}_0, v_0)\rightarrow 0$ almost surely, we simply confirm that the assumptions in Theorem~\ref{thm: BP SC} are met:  
Since $M=M_{0}=M_{n}$, Condition~\ref{C:Kuratowski} holds by Lemma~\ref{lem: kuratowski}. Since $(M,d)=(\mathbb{S}^{p-1},d)$ is compact, Condition~\ref{C:separable} and almost sure eventual boundedness of $\hat{S}^{0}$ is achieved. By Lemma~\ref{lem: continuous rho}, $d$ is equicontinuous. Also, $G(x)=x^{2}$ is $4$-subadditive, non-decreasing, and continuous. Thus, by Lemma~\ref{lem: locally integrable}, Conditions~\ref{C:rho-continuous} and \ref{C:locally-integrable} holds. Applying Theorem~\ref{thm: BP SC} and Lemma~\ref{lem: singletone}, we have $h_{d}(\{\hat{v}_0\}, \{v_0\})=d(\hat{v}_0, v_0)\rightarrow 0$ almost surely.

For $k\ge 1$, suppose that $h_{d}(\hat{S}^{i}, S^{i})\rightarrow 0$ almost surely for $i=0,\ldots, k-1$. Then,  $d(\hat{v}_{i}, v_{i})\rightarrow 0$ almost surely for $i=0,\ldots, k-1$ by Lemma~\ref{lem: hausdorff distance between great sphere converse}. In our framework, the geodesic subspheres $\hat{S}^{k}$ and $S^{k}$ are the sample and population generalized \fre means, respectively, with 
%
%
$$
T=\mathbb{S}^{p-1}, (M,d)=(\mathfrak{S}^{k},h_{d}), M_{0}=\mathfrak{S}^{k}(S^{k-1}), M_{n}=\mathfrak{S}^{k}(\hat{S}^{k-1}), \ \mbox{and} \ \cost(t,m)=d^{2}(t,m)
$$
for $t\in T, m\in M$. In this setting, we check the assumptions in Theorem~\ref{thm: BP SC}.

Since $d(\hat{v}_{i}, v_{i})\rightarrow 0$ almost surely for $i=0,\ldots, k-1$, Lemma~\ref{lem: compact domain} ensures that 
$M_{n}$ converges to $M_{0}$ in the sense of Kuratowski. Thus, Condition~\ref{C:Kuratowski} holds. By Lemma~\ref{lem: compact domain}, $\mathfrak{S}^{k}$ is compact. Thus, Condition~\ref{C:separable} holds and almost sure eventual boundedness of $\hat{E}_{n}$ is achieved. By Lemma~\ref{lem: continuous rho}, $d$ is equicontinuous (Definition~\ref{def: equicontinuous}). Also, $G(x)=x^{2}$ is $4$-subadditive, non-decreasing and continuous. Thus, by Lemma~\ref{lem: locally integrable}, Conditions~\ref{C:rho-continuous} and \ref{C:locally-integrable} hold. By Theorem~\ref{thm: BP SC} and Lemma~\ref{lem: singletone}, $h_{d}(\hat{S}^{k},S^{k})\rightarrow 0$ almost surely.

By induction principle, we  conclude that $h_{d}(\hat{S}^{k}, S^{k})\rightarrow 0$ for $k=0,\ldots, p-1$. As a byproduct,  Lemma~\ref{lem: hausdorff distance between great sphere converse} gives that $d(\hat{v}_{i}, v_{i})\rightarrow 0$ $(i=0,\ldots, p)$ almost surely as well.
\end{proof}


\subsection{Condition~\ref{C:Kuratowski} holds when \texorpdfstring{$M_n \equiv M_0$}{Mn = M0}}

Condition~\ref{C:Kuratowski}, a key requirement for BP consistency, is trivially satisfied when the minimizing domains are fixed, i.e., $M_n \equiv M_0$; see Lemma~\ref{lem: kuratowski} below. Thus, our result applies to various extensions of \fre means with fixed minimizing domains.

\begin{lemma}\label{lem: kuratowski}
If $M_{n}=M_{0} \in {CL}(M)$ for all $n \ge 1$, $M_n$ converges to $M_{0}$ in the sense of Kuratowski almost surely.
\end{lemma}

\begin{proof}[of Lemma~\ref{lem: kuratowski}]
Since $M_{0}$ is closed, for every sequence $m_{n}\in M_{n}=M_{0}$, accumulation point of $m_{n}$ is in $M_{0}$. Moreover, for every $m_{0}\in M_{0}$, sequence $m_{n}=m_{0} \in M_{n}=M_{0}$ converges to $m_{0}$. Thus, $M_{n}$ converges to $M_{0}$ in the sense of Kuratowski.  
\end{proof}

\subsection{Convergence of minima and existence of population generalized \fre mean}

In Theorem~\ref{thm: convergence of minimum} below, we establish two additional results: the convergence of empirical minima to the population minimum of the risk functions and the existence of population \fre means. This theorem is used in the proof of Corollary~\ref{cor: BP SC huckemann generalization}.

The minima of the objective functions $F_0$ and $F_n$ are given by 
$$l_{0}=\inf_{m\in M_{0}} F(m),\quad \hat{l}_{n}(\omega)=\inf_{m\in M_{n}} F_{n}(\omega)(m).$$ 

\begin{theorem}\label{thm: convergence of minimum}
Suppose Conditions~\ref{C:Kuratowski}, \ref{C:separable}, \ref{C:rho-continuous}, and \ref{C:locally-integrable} hold. Then almost surely, $$\limsup_{n\rightarrow \infty } \hat{l}_n(\omega) \le l_{0}.$$  
Additionally, assume Condition~\ref{C:Heine-Borel} and the following:

(i) $\hat{E}_{n}$ is eventually bounded almost surely, 

(ii) $\hat{E}_{n}$ is eventually non-empty almost surely, i.e., $$\mathcal{P}\left[\cup_{N=1}^{\infty}\cap_{n=N}^{\infty}\{\omega\in \Omega | \hat{E}_{n}(\omega) \neq  \emptyset\} \right]=1.$$

Then, $\lim_{n \to \infty} \hat{l}_n(\omega) = l_{0}$ almost surely and $E_{0}\neq \emptyset$.
\end{theorem}

\begin{proof}[of Theorem~\ref{thm: convergence of minimum}]
The first part, $\limsup_{n\rightarrow \infty } \hat{l}_n(\omega) \le l_{0}$ has been shown in the proof of Lemma~\ref{lem:Z-SC}. There, we established that, under Conditions~\ref{C:Kuratowski}, \ref{C:separable}, \ref{C:rho-continuous}, and \ref{C:locally-integrable},
$\limsup_{n\rightarrow \infty } \hat{l}_n(\omega) \le  l_{0}$ holds almost surely. 

To verify $\lim_{n \to \infty} \hat{l}_n(\omega) = l_{0}$ almost surely, it suffices to show that for almost every $\omega \in \Omega$, the following holds. For any subsequence $n_{a}$, there exists a further subsequence $n_{a(b)}$ such that $\hat{l}_{n_{a(b)}}(\omega)$ converges to $l_{0}$.

Fix $\omega \in \Omega$ satisfying the following conditions:
$M_n(\omega)$ converges to $M_0$ in the sense of Kuratowski, $\epsilon_{n}(\omega) \to 0$, $\hat{E}_{n}(\omega)$ is eventually bounded and eventually non-empty. The set of such $\omega$ has probability 1. We will show $\lim_{n\rightarrow \infty } \hat{l}_n(\omega) = l_{0}$ holds for such $\omega$.

Choose an arbitrary subsequence $n_{a}$, and let $m_{n_{a}} \in \hat{E}_{n_{a}}$. By the eventual boundedness of $\hat{E}_{n}$ and Condition~\ref{C:Heine-Borel}, there exists a subsequence of $n_a$, denoted by $n_{a(b)}$, such that $m_{n_{a(b)}}$ converges to some $m \in M$. By Theorem~\ref{thm: BP SC}, $m \in E_{0}$, thus $F(m) = l_{0}$. By Lemma~\ref{lem:epi-convergence}, $F_{n_{a(b)}}(\omega)(m_{n_{a(b)}})$ converges to $F(m) = l_{0}$. Since 
$$\hat{l}_{n_{a(b)}}(\omega) \leq F_{n_{a(b)}}(\omega)(m_{n_{a(b)}}) \leq \hat{l}_{n_{a(b)}}(\omega) + \epsilon_{n_{a(b)}}(\omega)$$
and $\epsilon_{n} \to 0$, it follows that $\hat{l}_{n_{a(b)}}(\omega) \to l_{0}$.  We have shown that $m\in E_{0}$, thus $E_{0}\neq \emptyset$.
\end{proof}


\subsection{Beyond Corollary~\ref{cor: BP SC Hd}: weaker moment assumptions and extension to \texorpdfstring{$p>0$}{p>0}}\label{app: further generalization} 

This subsection is motivated by the work of \cite{schotz2022strong}, which introduced shifted cost functions and lower-boundedness to relax moment assumptions and extend the scope of strong consistency results for $L^{p}$-center of mass. Building upon their framework, we extend Corollary~\ref{cor: BP SC Hd} in two directions: (1) weakening the moment condition from $E\{d^p(X,o)\}<\infty$ to $E\{d^{p-1}(X,o)\}<\infty$ for $p>1$, and (2) covering the regime $0<p\le 1$, which was not considered before in the strong consistency result of $L^{p}$-center of mass with random minimizing domian by \citep{evans2024limit}.

We begin by specifying a general class of cost functions that ensure the BP consistency of $\hat{E}_n$. In the additional result stated in Corollary~\ref{cor: BP SC schotz generalization} below, the cost functions are required to be \textit{lower-bounded}, a condition that may be overly technical to verify in typical applications.

\begin{definition}\label{def: lower bound} ~\citep{schotz2022strong}
We say that $\cost$ is lower-bounded if there exists $o \in M$, $\psi^{+}, \psi^{-}:[0, \infty) \rightarrow[0, \infty)$, $a^{+}, a^{-} \in(0, \infty)$, and $\sigma\left(X_1, \ldots, X_n\right)$-measurable random variables $a_n^{+}, a_n^{-} \in[0, \infty)$ such that
\begin{equation}
\label{eq:lowerboundSchotz}
\begin{aligned}
& a^{+} \psi^{+}\{d(m,o)\}-a^{-} \psi^{-}\{d(m,o)\} \leq F(m)=E\{\cost(X,m)\}, \\
& a_n^{+} \psi^{+}\{d(m,o)\}-a_n^{-} \psi^{-}\{d(m,o)\} \leq F_{n}(m) = \frac{1}{n}\sum_{i=1}^{n}  \cost(X_{i},m)
\end{aligned}
\end{equation} 
for all $m \in M$. Furthermore, ${a}_n^{+} \stackrel{n \rightarrow \infty}{\longrightarrow} {a}^{+}$ and ${a}_n^{-} \stackrel{n \rightarrow \infty}{\longrightarrow} {a}^{-}$ almost surely. Lastly, $\psi^{+}(\delta) / \psi^{-}(\delta) \rightarrow \infty$ as $\delta\rightarrow \infty$ and $\psi^{-}$ is non-decreasing.
\end{definition}

\begin{corollary}\label{cor: BP SC schotz generalization}
Assume Conditions~\ref{C:Kuratowski}, \ref{C:separable}, \ref{C:Heine-Borel}, \ref{C:rho-continuous}, and \ref{C:locally-integrable}. Additionally, assume that the cost function $\cost$ is lower-bounded in the sense of Definition \ref{def: lower bound}. Then, $\hat{E}_{n}$ is BP consistent with $E_{0}$. 
\end{corollary}

The lower-boundedness of cost functions was introduced in \cite{schotz2022strong} to establish the eventual boundedness of the empirical generalized \fre mean and its BP consistency for the case $M_{0} = M_n = M$. 
We refer to Remark 3.4 of \cite{schotz2022strong} for further explanation and examples of such cost functions.

\begin{remark} \textit{In Definition~\ref{def: lower bound}, we additionally assume that \emph{\(\psi^{-}\) is non-decreasing}, a condition not included the definition of lower-bounded functions in \cite{schotz2022strong}. We argue that this additional assumption is crucial for establishing the eventual boundedness of \(\hat{E}_{n}\). 
To be specific, let $\ell_0$ and $\ell_n$ denote the left-hand sides of (\ref{eq:lowerboundSchotz}). The eventual boundedness of \(\hat{E}_{n}\) relies on the divergence of  $F_n(m)$ as $m\in M$ moves away from a fixed point $o \in M$. Without the non-decreasing assumption on $\psi^{-}$,  $\ell_0$ and $\ell_n$ may fail to serve as effective lower bounds, which are necessary to ensure the divergence of $F_0$ and $F_n$. An example of functions $\psi^+$ and $\psi^-$ satisfying all requirements of Definition~\ref{def: lower bound}  except for the non-decreasing property of $\psi^-$ is as follows: 
\begin{align*}
    \psi^{+}(\delta) & = I(\delta < 1) + \frac{1}{\delta} I(\delta \ge 1), \\ 
     \psi^{-}(\delta) &= I(\delta < 1) + \frac{1}{\delta^{2}} I(\delta \ge 1).
\end{align*}
In this case, $\ell_0$ and $\ell_n$ do not diverge  as $d(m,o)\rightarrow \infty$, making them unsuitable for verifying the eventual boundedness of 
$\hat{E}_{n}$.  
All proofs in \cite{schotz2022strong} that rely on lower-boundedness remain valid with this additional assumption on $\psi^{-}$. Furthermore, the examples of lower-bounded cost functions given in \cite{schotz2022strong} also satisfy this non-decreasing criterion for \(\psi^{-}\).}



\end{remark}

Corollary~\ref{cor: BP SC schotz generalization} extends Theorem 3.5 of \cite{schotz2022strong} to cases where the minimizing domain is random. Additionally, our result relaxes certain conditions imposed in \cite{schotz2022strong}. A detailed comparison between Theorem 3.5 of \cite{schotz2022strong} and Corollary~\ref{cor: BP SC schotz generalization} is provided in Section~\ref{SM: comparison with previous ones}.



We now provide a proof of Corollary~\ref{cor: BP SC schotz generalization}. To establish BP consistency of $\hat{E}_{n}$ in this setting, it suffices to show that the assumption of lower-bounded cost functions ensures the almost sure eventual boundedness of $\hat{E}_{n}$ (\textit{cf}. Theorem \ref{thm: BP SC}). 


\begin{lemma}\label{lem: lower bound}
Assume Conditions~\ref{C:Kuratowski}, \ref{C:separable}, \ref{C:rho-continuous}, and \ref{C:locally-integrable}. If $\cost$ is lower-bounded, then $\hat{E}_{n}$ is eventually bounded almost surely.
\end{lemma}

\begin{proof}[of Lemma~\ref{lem: lower bound}]
By Theorem~\ref{thm: convergence of minimum}, ${\limsup}_{{n\rightarrow\infty}}\hat{l}_{n}(\omega)\le l_{0}$ almost surely. $a_{n}^{+}\rightarrow a^{+}$ and $a_{n}^{-}\rightarrow a^{-}$ almost surely by assumption, where $a_{n}^{+},a^{+},a_{n}^{-}$, and $a^{-}$ are defined in Definition~\ref{def: lower bound}. Fix $\omega\in \Omega$ satisfying ${\limsup}_{{n\rightarrow\infty}}\hat{l}_{n}(\omega)\le l_{0}$, $a_{n}^{+}\rightarrow a^{+}$, $a_{n}^{-}\rightarrow a^{-}$, and $\epsilon_{n}(\omega) \rightarrow 0$. The set of such $\omega$ has probability $1$.

We will show that $\hat{E}_{n}(\omega)$ is eventually bounded for such $\omega$. For such purpose, it is enough to show that there exists a bounded subset $B(\omega)$ of $M$ and 
$N(\omega)\ge 1$ such that $\inf _{m \in M_{n}(\omega)\backslash B(\omega)}  F_n(m)(\omega)  > \hat{l}_{n}(\omega)+\epsilon_{n}$ for $n> N(\omega)$. This implies $\hat{E}_{n}(\omega) \subset B(\omega)$ for all $n> N(\omega)$, thus $\hat{E}_{n}(\omega)$ is eventually bounded.

We fix $0<\epsilon< a^{+}$. Since $\lim_{\delta\rightarrow\infty} \psi^{+}(\delta)/\psi^{-}(\delta) = \infty$ and $\psi^{-}$ is non-decreasing, we have
\begin{equation*}
\begin{split}
& \lim_{\delta\rightarrow\infty} ({a}^{+}-\epsilon)\psi^{+}(\delta) - ({a}^{-}+\epsilon)\psi^{-}(\delta) = \lim_{\delta\rightarrow\infty} \psi^{-}(\delta) \{ ({a}^{+}-\epsilon)\psi^{+}(\delta)/\psi^{+}(\delta) - ({a}^{-}+\epsilon) \} = \infty.
\end{split}    
\end{equation*}
Thus, there exist $\delta_{o}>0$ such that $\inf_{\delta \ge \delta_{o}} ({a}^{+}-\epsilon) \psi^{+}(\delta)-({a}^{-}+\epsilon) \psi^{-}(\delta) \ge l_{0}+2$. We set $B(\omega)=B(o, \delta_{o})$, for $o \in M$ in Definition~\ref{def: lower bound}. 

Choose $N_{1}(\omega)$ such that $\mid a_{n}^{+}-a^{+} \mid \le \epsilon$ and $\mid a_{n}^{-}-a^{-} \mid \le \epsilon$ for $n > N_{1}(\omega)$. 
Since 
\begin{align*}
    \inf _{m \in M_{n}(\omega)\backslash B(\omega)} F_{n}(\omega)(m) & \ge \inf_{\delta\ge \delta_{0}} a_{n}^{+}\psi^{+}(\delta)-a_n^{-} \psi^{-}(\delta) \\
    &\ge \inf_{\delta\ge \delta_{o}}(a^{+}-\epsilon)\psi^{+}(\delta)-(a^{-}+\epsilon) \psi^{-}(\delta) \\
    &\ge l_{0}+2
\end{align*} 
for all $n > N_{1}(\omega)$,
we have $\inf _{m \in M\backslash B(\omega)} F_{n}(\omega)(m) \ge l_{0}+2$ for $n > N_{1}(\omega)$.

Since $\limsup_{n\rightarrow \infty} \hat{l}_{n}(\omega)\le l_{0}$, $\hat{l}_{n}(\omega)\le l_{0}+1$ eventually. Thus, we can choose $N_{2}(\omega)$ satisfying $\hat{l}_{n}(\omega)\le l_{0}+1$ for all $n > N_{2}(\omega)$. Since $\epsilon_{n}(\omega) \rightarrow 0$ as $n\rightarrow \infty$, we can choose $N_{3}(\omega)$ such that $\epsilon_{n}(\omega)<1$ for $n>N_{3}(\omega)$. Then for $n > N(\omega)=\max\{N_{1}(\omega),N_{2}(\omega),N_{3}(\omega)\}$, $\inf _{m \in M_{n}(\omega)\backslash B(\omega)}  F_n(m)(\omega) \geq l_{0}+2 > \hat{l}_{n}(\omega)+\epsilon_{n}$. 
\end{proof}

Corollary~\ref{cor: BP SC schotz generalization} follows from a combination of Theorem~\ref{thm: BP SC} and Lemma \ref{lem: lower bound}.

\begin{proof}[of Corollary~\ref{cor: BP SC schotz generalization}]
By Lemma~\ref{lem: lower bound}, $\hat{E}_{n}$ is eventually bounded almost surely. By Theorem~\ref{thm: BP SC}, $\hat{E}_{n}$ is BP consistent with $E_{0}$. 
\end{proof}

We now introduce a class of lower-bounded cost functions constructed by integrating a non-decreasing, $b$-subadditive function \(h\). Specifically, we consider cost functions of the form \(\cost(t, m) = H\{ d(t, m)\} - H \{ d(t, o) \}\), where \(H(t) = \int_0^t h(x)\, dx\). The following corollary ensures BP consistency under this general setting.

\begin{corollary}\label{cor: H-fre mean}
Suppose   \(T = M\), and  that the cost function is defined as \(\cost(t, m) = H\{ d(t, m)\} - H \{ d(t, o) \} \), where $o \in M$ and
$H$ is given by \(H(t) = \int_{0}^{t} h(x) dx\) ($t\ge 0$). Here, \(h\) is a non-decreasing, {$b$-subadditive} function that satisfies $\lim_{x \rightarrow \infty} h(x) = \infty$, and $E[h\{ d(X,o)\}]<\infty$. Furthermore, if Conditions~\ref{C:Kuratowski}, \ref{C:separable}, and \ref{C:Heine-Borel} hold, then \(\hat{E}_{n}\) is BP consistent with \(E_{0}\).
\end{corollary}

Recall that non-decreasing $h:[0,\infty)\rightarrow [0,\infty)$ is $b$-subadditive when there exists $b\ge 1$ such that $h(2x)\le b h(x)$ for all $x\ge 0$.

Corollary~\ref{cor: H-fre mean} can be shown by applying Corollary~\ref{cor: BP SC schotz generalization} to the cost function $\cost(t,m)=H\{d(t,m)\}-H\{d(t,o)\}$, by showing that it is lower-bounded and satisfies Conditions~\ref{C:rho-continuous} and \ref{C:locally-integrable}.

To this end, we collect basic properties of functions of the form $H(x)=\int_{0}^{x} h(y) dy$ in the following lemma. While these properties were previously established in \cite{schotz2022strong}, we provide proof here for completeness of the paper. 


\begin{lemma}[Lemma A.7 of \cite{schotz2022strong}]\label{lem: H properties}
Let $H: [0,\infty) \rightarrow [0,\infty) $ be defined by $H(x)=\int_{0}^{x} h(y) dy$, 
where  $h : [0,\infty) \rightarrow [0,\infty)$ is non-decreasing and $b$-subadditive. Then, the following properties hold.

(i) For $x,y\ge 0$, $\frac{1}{2}h(x)+\frac{1}{2}h(y) \le h(x+y)\le bh(x)+bh(y)$, where $b$ is the coefficient used in the definition of $b$-subadditivity. 

(ii) For $x\ge 0, H(x) \ge (x/2)h(x/2)$. 

(iii) $H(\mid x-y \mid)-H(x)\ge b^{-1} H(y)- 2yh(x)$.
\end{lemma}

\begin{proof}[of Lemma~\ref{lem: H properties}]
The first and second parts are immediate as follows. For $x,y\ge 0$, $$h(x)+h(y) \le 2 \max\{h(x),h(y)\} \le 2 h(x+y)$$ and $$h(x+y)\le h\{\max(2x,2y)\} \le h(2x)+h(2y) \le bh(x)+bh(y)$$ holds. 

For $x\ge 0$, $$H(x)=\int_{0}^{x} h(t) dt \ge \int_{x/2}^{x} h(t) dt \ge h(x/2)(x/2).$$    

For the third part, we first consider the case \( x \ge y \). Define \( f(x, y) = H(x - y) - H(x) - b^{-1}H(y) + 2y h(x) \). We want to show \( f(x, y) \ge 0 \). The derivative of \( f \) with respect to \( y \) is
\[
\partial_y f(x, y) = -h(x - y) - b^{-1}h(y) + 2h(x).
\]
Since $2h(x)\ge h(x-y)+h(y)$ and $b^{-1}\le 1$, we obtain \( \partial_y f(x, y) \ge 0 \) as \( b^{-1} \le 1 \). Hence, \( f(x, y) \ge f(x, 0) = 0 \), as \( H(0) = 0 \).

Now, consider the case \( x \le y \). Set \( g(x, y) = H(y - x) - H(x) - b^{-1}H(y) + 2yh(x) \), which yields
\[
\partial_y g(x, y) = h(y - x) - b^{-1}h(y) + 2h(x).
\]
By substituting $h(y) \le bh(y-x)+bh(x)$, we obtain \( \partial_y g(x, y) \ge h(x) \ge 0 \). Thus, \( g(x, y) \ge g(x, x) = -(1 + b^{-1})H(x) + 2xh(x) \) as \( H(0) = 0 \). By the definition of \( H \), as \( h \) is non-decreasing, \( H(x) \le xh(x) \). Hence, \( g(x, y) \ge 0 \) as \( 1 + b^{-1} \le 2 \).

Together, we have shown \( H(|x - y|) - H(x) - b^{-1}H(y) + 2yh(x) \ge 0 \) for all \( x, y \ge 0 \).
\end{proof}

Next intermediate result in Lemma~\ref{lem: H-fre mean} builds upon Lemma~\ref{lem: H properties}, and establishes that cost functions of the form $\cost(t,m) = H\{d(t,m)\}-H\{d(t,o)\}$ satisfy Conditions~\ref{C:rho-continuous} and \ref{C:locally-integrable}, and are lower-bounded. 

\begin{lemma}[Corollaries 4.3 and 4.4 of \cite{schotz2022strong}]\label{lem: H-fre mean}
Suppose $T=M$ and that the cost function can be written as $\cost(t,m) = H\{d(t,m)\}-H\{d(t,o)\}$ for a fixed $o\in M$, where $H: [0,\infty) \rightarrow [0,\infty)$ is given  by $H(x)=\int_{0}^{x} h(y) dy$ with non-decreasing and $b$-subadditive $h : [0,\infty) \rightarrow [0,\infty)$ which satisfies $\lim_{x\rightarrow \infty} h(x) = \infty$ and $E [h\{d(X,o)\}]<\infty$. Then $E \{ \cost(X,m) \}<\infty $ for all $m\in M$, Conditions~\ref{C:rho-continuous} and \ref{C:locally-integrable} are satisfied, and $\cost(t,m)$ is lower-bounded.
\end{lemma}


\begin{proof}[of Lemma~\ref{lem: H-fre mean}]
We first show Condition~\ref{C:rho-continuous}. Since $H$ is defined by integrating $h$, $H$ is continuous. Moreover, $d(t,m)$ is a continuous function of $m$ for each fixed $t\in T$. Thus $\cost(t, \cdot)=H\{d(t,\cdot)\}-H\{d(t,o)\}: M\rightarrow \mathbb{R}$ is continuous for each $t\in T$.

Secondly, we show that $\cost(X,m)=H\{d(X,m)\}-H\{d(X,o)\}$ is integrable for each $m\in M$, when $h\{d(X,o)\}$ is integrable. By applying Lemma~\ref{lem: H properties}, we have
\begin{align*} 
| H\{d(X,m)\}-H\{ d(X,o)\} | & \le | d(X,m)-d(X,o) | \times  h[  \max\{d(X,m), d(X,o)\}] \\
& \le d(m,o) h\{d(m,o)+d(X,o)\} \\
& \le b d(m,o)[h\{d(m,o)\}+ h\{d(X,o)\}].  
\end{align*} 
Thus, $\cost(X,m)=H\{ d(X,m) \}-H \{ d(X,o) \}$ is integrable for each $m\in M$.

Third, we show Condition~\ref{C:locally-integrable}. Since $\cost(X,m')-\cost(X,m) = H\{ d(X,m') \}-H\{d(X,m)\}$, we get 
\begin{align*} 
 |\cost(X,m')-\cost(X,m) | &\le | d(X,m')-d(X,m) | \times h [\max\{d(X,m'), d(X,m)\}] \\
& \le d(m,m') h \{ d(m',m)+d(X,m) \} \\
&\le bd(m,m') [ h\{ d(m',m) \}+ h\{ d(X,m) \}], 
\end{align*}
using a similar argument as above. Therefore, $\pi_{m,r}(X)$ and $\Pi_{m,r}(X)$ are integrable for all $r>0$.

Finally, the lower-boundedness of $\cost$ is now established. By applying Lemma~\ref{lem: H properties}, we obtain
\begin{align*} 
 H\{d(t,m)\}-H\{d(t,o)\} & \ge  H\{\mid d(t,o)-d(o,m) \mid\}-H\{d(t,o)\} \\
& \ge b^{-1}H\{d(o,m)\}-2h\{d(t,o)\}d(o,m) \\
&\ge b^{-1} h\{d(o,m)/2\} d(o,m)/2-2h\{d(t,o)\} d(o,m). 
\end{align*}
Thus, by setting $\psi^{+}(\delta)=b^{-1}h(\delta/2)(\delta/2)$, $\psi^{-}(\delta)=2\delta$, $a^{+}=a_{n}^{+}=1$, $a^{-}=E [ h\{d(X,o)\}]$, $a_{n}^{-}=\frac{1}{n} \sum_{i=1}^{n} h\{ d(X_{i},o)\}$, we can easily check that $\cost$ is lower-bounded.    
\end{proof}

Combining the results of Corollary~\ref{cor: BP SC schotz generalization} and Lemma~\ref{lem: H-fre mean}, we directly obtain Corollary~\ref{cor: H-fre mean}.

\begin{proof}[of Corollary~\ref{cor: H-fre mean}]
By Lemma~\ref{lem: H-fre mean}, Conditions~\ref{C:rho-continuous} and \ref{C:locally-integrable} holds and $\hat{E}_{n}$ is eventually bounded almost surely. By Corollary~\ref{cor: BP SC schotz generalization}, $\hat{E}_{n}$ is BP consistent with $E_{0}$.
\end{proof}

Finally, we generalize to the case where $\cost(t,m) = d^p(t, m) - d^{p}(t, o)$.

\begin{corollary}\label{cor: alpha-fre mean combined}
Suppose $T = M$, and for $p > 0$, $\cost(t,m) = d^p(t, m) - d^{p}(t, o)$, where $o \in M$. For the case $p > 1$, assume in addition that $d^{p - 1}(X, o)$ is integrable. If  Conditions~\ref{C:Kuratowski}, \ref{C:separable}, and \ref{C:Heine-Borel} hold, and $\hat{E}_{n}$ is eventually bounded almost surely, then \(\hat{E}_{n}\) is BP consistent with \(E_{0}\).
\end{corollary}

\begin{proof}[of Corollary~\ref{cor: alpha-fre mean combined}] When \( p > 1 \), it is straightforward to verify that \( h(x) = p x^{p - 1} \) is non-decreasing, $\lim_{x\rightarrow \infty} h(x) = \infty$, and $b$-subadditive with \( b = 2^{p - 1} > 1 \). Thus, the conclusion of Corollary~\ref{cor: alpha-fre mean combined} directly follows from Corollary~\ref{cor: H-fre mean} when $p>1$.

For the case \( p \le 1 \), we can find a function \( h : [0, \infty) \rightarrow [0, \infty) \) that is strictly increasing (and therefore non-decreasing), continuous, concave (implying $b$-subadditivity with \( b = 2 \)), $\lim_{x\rightarrow \infty} h(x) = \infty$, and satisfies \( E [h\{ d^{p}(X, o)\}] < \infty \) \citep[Lemma A.8]{schotz2022strong}. For such $h$, \( H(x) = \int_{0}^{x} h(t) \, dt \) is convex and strictly increasing, and thus has an inverse \( H^{-1} \) that is strictly increasing, $H^{-1}(0)=0$ and concave. We define a new metric \( d_{H} = H^{-1}(d^{p}) \) on \( M \). Then our cost function \( \cost(t, m) = d^{p}(t, m)-d^{p}(t, o) \) can be represented by \( H\{d_{H}(t, m)\} - H\{d_{H}(t, o)\}\). 

Therefore, to prove the result in Corollary~\ref{cor: alpha-fre mean combined} for \( p \le 1 \) using Corollary~\ref{cor: H-fre mean}, we need to verify \( E [h\{d_{H}(X, o)\}] < \infty \) and Conditions~\ref{C:Kuratowski}--\ref{C:locally-integrable} still hold and $\hat{E}_{n}$ is almost surely eventually bounded for the new metric space \( (M, d_{H})\). Then by Corollary~\ref{cor: H-fre mean}, $\hat{E}_{n}$ is BP consistent with $E_{0}$, with respect to the new metric $d_{H}$. Since $d_{H}(m,m')=H^{-1}\{d^{p}(m,m')\}$, we have
$$\underset{n\rightarrow \infty}\lim\sup_{m\in \hat{E}_{n}} d_{H}(m,E)=0$$ if and only if $$\underset{n\rightarrow \infty}\lim\sup_{m\in \hat{E}_{n}} d(m,E)=0.$$ Thus, $\hat{E}_{n}$ is also BP consistent with $E_{0}$ with respect to $d$, we obtain the conclusion of Corollary~\ref{cor: alpha-fre mean combined}.

We first verify that \( (M, d_{H})\) is a metric space. It is trivial that $d_{H}(m,m')\ge 0$ and $d_{H}(m,m')=0$ if and only if $d_{H}(m,m')=0$ for $m,m'\in M$, $d_{H}(m,m')=d_{H}(m',m)$. Moreover, since $H^{-1}$ and $x\mapsto x^{p} : [0,\infty) \rightarrow [0,\infty)$ is strictly increasing and concave, we get 
\begin{align*} 
 d_{H}(m,m'') & =H^{-1}\{d^{p}(m,m'')\} \\
&\le H^{-1}[\{d(m,m')+d(m',m'')\}^{p}]\\
& \le H^{-1}\{d^{p}(m,m')+d^{p}(m',m'')\}\\
& \le H^{-1}\{d^{p}(m,m')\}+H^{-1}\{d^{p}(m',m'')\}\\
& = d_{H}(m,m')+d_{H}(m',m'')     
\end{align*}
for all $m,m',m''\in M$. Thus, \( (M, d_{H})\) is a metric space.  

Since $H^{-1}$ is concave, we get \( H^{-1}(x) \le Ax+B \) for some \( A,B \ge 0 \). Since $h$ is concave and non-decreasing, we have $$ h(Ax + B) \le h(Ax)+h(B) \le \max\{A,1\} h(x)+h(B). $$ 
Thus,  
\begin{equation*}
\begin{split}
E [h\{d_{H}(X, o)\}] & = E [h\circ H^{-1}\{d^{p}(X, o)\}]\\
& \le E [h\{Ad^{p}(X,o)+B\}]  \\
& \le AE[h\{d^{p}(X,o)\}]+h(B) < \infty.
\end{split}    
\end{equation*} 

To verify Conditions~\ref{C:Kuratowski}--\ref{C:locally-integrable} and almost sure eventual boundedness of $\hat{E}_{n}$ still hold for new metric space $(M,d_{H})$, we show that $H^{-1}: [0,\infty) \rightarrow [0,\infty)$ is a one to one correspondence between $[0,\infty)$ and $[0,\infty)$.

We can observe that $H^{-1}(0) = 0$, $H^{-1}$ is strictly increasing, continuous, and satisfies $\lim_{x \to \infty} H^{-1}(x) = \infty$, $H^{-1}: [0,\infty) \to [0,\infty)$. Therefore, $H^{-1}: [0,\infty) \rightarrow [0,\infty)$ is a one to one correspondence between $[0,\infty)$ and $[0,\infty)$. As a result, the metric topology induced by $(M, d)$ and the metric topology induced by $(M, d_{H})$ are equivalent. Moreover, subset $C\subset M$ is bounded in $(M, d)$ if and only if it is bounded in $(M, d_H)$. Thus, Conditions~\ref{C:Kuratowski}--\ref{C:locally-integrable}, and the almost sure eventual boundedness of $\hat{E}_{n}$ are preserved.
\end{proof}

\section{Comparison with previous results on the law of large numbers for various extensions of \fre means}
\label{SM: comparison with previous ones}

\subsection*{Outline}

In this section, we review previous results on the strong consistency of various extensions of \fre means, including \fre $\rho$-means~\citep{huckemann2011intrinsic}, $\cost$-\fre means~\citep{schotz2022strong}, and $C$-restricted \fre $p$-means~\citep{evans2024limit}. We then explore their relationship with our Corollaries~\ref{cor: BP SC huckemann generalization}, \ref{cor: BP SC schotz generalization}, \ref{cor: H-fre mean}, and \ref{cor: alpha-fre mean combined}.

\subsection{\texorpdfstring{\fre $\rho$-mean of \cite{huckemann2011intrinsic}}{Frechet rho-mean of Huckemann et al. (2011)}}

We review the strong consistency results for the \fre $\rho$-mean established in \cite{huckemann2011intrinsic}, and compare them with Corollary~\ref{cor: BP SC huckemann generalization} in this paper. As discussed in Section~\ref{seC:frechet means}, the \fre $\rho$-mean is a special case of our generalized \fre mean, obtained by setting $\cost = \rho^2$ for some function $\rho: T \times M \to [0, \infty)$ and assuming fixed minimizing domains, i.e., $M = M_0 = M_n$ for all $n \geq 1$ and $\varepsilon_n = 0$.

In \cite{huckemann2011intrinsic}, strong consistency is shown in two steps: Ziezold consistency under the equicontinuity of $\rho$, and Bhattacharya-Patrangenaru (BP) consistency under its coercivity. We summarize their results below within our framework.

\begin{theorem}[\fre $\rho$-Mean \citep{huckemann2011intrinsic}]\label{cor: BP SC huckemann} Let $\cost = \rho^2$ with $\rho: T \times M \to [0, \infty)$, $M = M_{0} = M_n$ for all $n \geq 1$ and $\omega \in \Omega$, and $\varepsilon_n = 0$ for all $n \geq 1$. 
Assume $(M, d)$ is separable, $E \{\rho^{2}(X, m)\} < \infty$ for all $m \in M$, $\rho$ is continuous on $T \times M$, equicontinuous and coercive, $E_{0} \neq \emptyset$, $\cup_{n=1}^{\infty} E_{n}(\omega)$ satisfies the Heine-Borel property almost surely. Then, $\hat{E}_n$ is BP consistent with $E_{0}$.
\end{theorem}

We now discuss the common and differing assumptions between Theorem~\ref{cor: BP SC huckemann} and our result (Corollary~\ref{cor: BP SC huckemann generalization}). The primary difference is that Corollary~\ref{cor: BP SC huckemann generalization} allows for random minimizing domains $M_n$, whereas Theorem~\ref{cor: BP SC huckemann} assumes a fixed domain $M$. This generalization expands the range of applicable scenarios.

Even in the fixed-domain setting ($M_n \equiv M_0$), our corollary accommodates a larger class of cost functions. Specifically, it allows cost functions of the form $\cost = G \circ \rho$, where $G: [0, \infty) \to [0, \infty)$ is non-decreasing, continuous, $b$-subadditive, and satisfies $\lim_{x\to\infty} G(x) = \infty$. The classical choice $G(x) = x^2$ used in \cite{huckemann2011intrinsic} is just one instance of this more general setting.

Both results share standard assumptions, such as the separability of $(M, d)$ and the (almost sure) Heine-Borel property of the minimizing sets. However, the random domain setting in our result introduces a few additional technical requirements. For instance, Corollary~\ref{cor: BP SC huckemann generalization} assumes completeness of $(M, d)$, which is crucial for handling the measurability of bounding functions (see Lemma~\ref{lem: scalarization properties}), and for generalizing the cost beyond squared distance.

Our framework does not require continuity of the cost function in the product topology of $  T\times M$, which is needed in Theorem~\ref{cor: BP SC huckemann}. Instead, it relies on Conditions~\ref{C:Kuratowski}--\ref{C:locally-integrable}, including Kuratowski convergence of the random domains and local integrability of $\pi_{m,r}(X)$, $\Pi_{m,r}(X)$.

Finally, while Corollary~\ref{cor: BP SC huckemann generalization} assumes the Heine-Borel property of the entire space $M$, this can be weakened to the almost sure eventual boundedness of $\hat{E}_n$ together with the Heine-Borel property of their supports, as is standard in BP consistency arguments.

In summary, Corollary~\ref{cor: BP SC huckemann generalization} generalizes the strong consistency result of \cite{huckemann2011intrinsic} in two major directions: it permits random minimizing domains and supports a broader class of cost functions. Although this requires a few additional assumptions, these are necessary to ensure measurability and convergence under the more flexible setting.

\subsection{\texorpdfstring{$\mathfrak{c}$}{\textit{c}}-Fr\'{e}chet mean of \cite{schotz2022strong}}

As discussed in Section~\ref{seC:frechet means}, the $\cost$-\fre mean is a specific case of the generalized \fre mean, where \(M = M_{0} = M_n\) for all \(n \geq 1\). The work of \cite{schotz2022strong} established that $\hat{E}_n$ is a Ziezold and BP consistent estimator of \(E_{0}\) under conditions that share similarities with our framework. To maintain completeness, we restate these results within our setting.

\begin{theorem}[Theorems 3.2 and 3.5 of \cite{schotz2022strong}]
\label{thm: Z SC, BP SC schotz}
Suppose \(M_{0} = M_n = M\) for all \(n \geq 1\) and \(\omega \in \Omega\), and that Condition~\ref{C:separable} holds. Additionally, assume that $\cost(t, \cdot): M \to \mathbb{R}$ is lower semi-continuous for all \(t \in T\), and that \(\mathbb{E} \{\inf_{m \in M} \cost(X, m)\} > -\infty\). Then, $\hat{E}_n$ is Ziezold consistent with \(E_{0}\).

Furthermore, assume $\cup_{n\ge 1} \hat{E}_n$ satisfies the Heine-Borel property almost surely, \(E \{\sup_{m \in B}|\mathfrak{c}(X, m)|\}<\infty\) for all bounded sets \(B\), and $\cost$ is lower-bounded (Definition~\ref{def: lower bound}). Then $\hat{E}_n$ is BP consistent with \(E_{0}\).
\end{theorem}

We now compare this result with our Corollary~\ref{cor: BP SC schotz generalization}, focusing on the similarities, differences, and implications of each.

Both results assume that \((M, d)\) is Polish (Condition~\ref{C:separable}) and $\cost$ is lower-bounded. The Heine-Borel property is crucial in establishing BP consistency in both cases.

The key differences lie in the conditions imposed on the cost function and the flexibility of the domain. Theorem~\ref{thm: Z SC, BP SC schotz} assumes that $\cost(t, \cdot)$ is only lower semi-continuous, whereas our Corollary~\ref{cor: BP SC schotz generalization} requires both upper and lower semi-continuity. On the other hand, while Corollary~\ref{cor: BP SC schotz generalization} requires local infinium and local supremum of $\cost$ are integrable, Theorem~\ref{thm: Z SC, BP SC schotz} assumes \(E \{\inf_{m \in M} \cost(X, m)\} > -\infty\) and $E \{\sup_{m \in B}|\mathfrak{c}(X, m)|\}<\infty$ for all bounded sets \(B\), which is strictly stronger than Condition~\ref{C:locally-integrable}. Finally, Corollary~\ref{cor: BP SC schotz generalization} generalizes the fixed-domain assumption of Theorem~\ref{thm: Z SC, BP SC schotz} by allowing the minimizing domain \(M_n\) to vary with \(n\), enabling the analysis of more complex or data-driven model classes.

In summary, Corollary~\ref{cor: BP SC schotz generalization} extends the results of \cite{schotz2022strong} by incorporating random minimizing domains and relaxing the integrability assumptions on the cost function. The requirement of upper semi-continuity in our framework is a necessary trade-off that enables this generalization.

\subsection{\texorpdfstring{$H$-\fre mean and \fre$\text{-}p$ means}{H-Fre mean and Fre-p means} of \cite{schotz2022strong}}

In addition to their general consistency result (Theorem~\ref{thm: Z SC, BP SC schotz}), \cite{schotz2022strong} also provided concrete applications to specific classes of cost functions, namely the $H$-\fre mean and the \fre-$p$ mean. These results are presented as Corollaries 4.3, 4.4, 5.1, and 5.2 in their work. We restate them here using the terminology of our generalized \fre mean framework.

\begin{corollary}[Corollaries 4.3 and 4.4 of \cite{schotz2022strong}]\label{cor: H-fre mean schotz}
Let $T=M$, $M=M_{0}=M_{n}$ for all $n\ge 1$, and $\cost(t,m) = H\{d(t,m)\}-H\{d(t,o)\}$ for a fixed $o\in M$, where $H: [0,\infty) \rightarrow [0,\infty)$ is defined by $H(t)=\int_{0}^{t} h(x) dx$ for non-decreasing, $\lim_{x\rightarrow \infty } h(x) = \infty$, $b$-subadditive $h : [0,\infty) \rightarrow [0,\infty)$ that satisfies $E [h\{ d(X,o) \}]<\infty $. Assume Condition~\ref{C:separable}. Then, the following holds.

$\hat{E}_{n}$ is Ziezold consistent with $E_{0}$. Additionally, if $\cup_{n\ge 1} \hat{E}_{n}$ satisfies the sample Heine-Borel property almost surely, then $\hat{E}_{n}$ is BP consistent with $E_{0}$.
\end{corollary}

\begin{corollary}[Corollaries 5.1 and 5.2 of \cite{schotz2022strong}]\label{cor: alpha-fre mean combined schotz}
Let \( T = M \), $M=M_{0}=M_{n}$ and \( \cost(t, m) = d^{p}(t, m) - d^{p}(t, o) \) for a fixed \( o \in M \) and \( p > 0 \). Assume Condition~\ref{C:separable} holds. Additionally, if \( p > 1 \), assume that \( d^{p - 1}(X, m) \) is integrable for some \( m \). Then, the following holds. 

\( \hat{E}_{n} \) is Ziezold consistent with \( E_{0} \). If $\cup_{n\ge 1} \hat{E}_{n}$ satisfies the sample Heine-Borel property almost surely, then \( \hat{E}_{n} \) is BP consistent with \( E_{0} \).
\end{corollary}

The only difference between Corollaries~\ref{cor: H-fre mean}, \ref{cor: alpha-fre mean combined} and Corollaries~\ref{cor: H-fre mean schotz}, \ref{cor: alpha-fre mean combined schotz} lies in the adaptation of the random minimizing domain, while all other conditions remain the same. Therefore, our Corollary~\ref{cor: H-fre mean} and Corollary~\ref{cor: alpha-fre mean combined} serve as generalizations of Corollaries 4.3, 4.4, 5.1, and 5.2 from \cite{schotz2022strong}.

\subsection{C-restricted {Fr\'{e}chet} p-mean of \texorpdfstring{\cite{evans2024limit}}{Evans et al. (2024)} }

As discussed in Section~\ref{seC:frechet means}, the $C$-restricted \fre mean of \cite{evans2024limit} is a special case of our generalized \fre mean, when $T=M$ and $\cost(t,m)=d^{p}(t,m)$ for $p\ge 1$. \cite{evans2024limit} have shown that empirical $C_{n}$-restricted \fre $p$-mean is a Ziezold and BP consistent estimator of population $C$-restricted \fre $p$-mean under regularity conditions. For reference, we state their result below. 

\begin{theorem}[Theorems 4.3 and 4.4 of \cite{evans2024limit}]\label{thm: C-restricted fre means} Let $T=M$ and $\cost(t,m)=d^{p}(t,m)$ for $p\ge 1$. Assume Condition~\ref{C:Kuratowski} and $(M,d)$ is a separable metric space. If $E \{ d^{p-1}(X,m)\}<\infty$ for all $m\in M$, then $\hat{E}_{n}$ is a Ziezold consistent estimator of $E_{0}$. If $E \{ d^{p}(X,m) \}<\infty$ for all $m\in M$ and $(M,d)$ satisfies Heine-Borel property, then $\hat{E}_{n}$ is BP consistent with $E_{0}$.
\end{theorem}

We now compare the special case of our framework, stated in Corollary~\ref{cor: alpha-fre mean combined}, with the main result of \cite{evans2024limit} (Theorem~\ref{thm: C-restricted fre means}). 
Both results assume Condition~\ref{C:Kuratowski} and the separability of \( (M, d) \), establishing a common foundation for their respective consistency results.

The key advantage of our result is its weaker integrability condition and applicability to a broader range of $p>0$. \cite{evans2024limit} (Theorem~\ref{thm: C-restricted fre means}) requires \( E\{ d^{p-1}(X, m)\} < \infty \) for Ziezold consistency and \( E\{ d^{p}(X, m) \} < \infty \), along with the Heine-Borel property, for BP consistency, where \( p > 1 \). In contrast, our result (Corollary~\ref{cor: alpha-fre mean combined}) imposes the integrability of \( d^{p - 1}(X, m) \) for (\( p > 1 \)), making it a more relaxed condition than in \cite{evans2024limit}. Moreover, Corollary~\ref{cor: alpha-fre mean combined} extends the applicability to a broader range, allowing for cases where $0<p\le 1$.

\bibliographystyle{asa}
\bibliography{library}

@Article{Bhattacharya2003,
  author    = {Bhattacharya, Rabi and Patrangenaru, Vic},
  journal   = {The Annals of Statistics},
  title     = {{Large Sample Theory of Intrinsic and Extrinsic Sample Means on Manifolds. I}},
  year      = {2003},
  issn      = {00905364},
  month     = feb,
  number    = {1},
  pages     = {1--29},
  volume    = {31},
  publisher = {Institute of Mathematical Statistics},
  url       = {http://www.jstor.org/stable/3448366},
}

@Article{damon2014backwards,
  author    = {Damon, James and Marron, JS},
  title     = {Backwards principal component analysis and principal nested relations},
  journal   = {Journal of Mathematical Imaging and Vision},
  year      = {2014},
  volume    = {50},
  number    = {1-2},
  pages     = {107--114},
  publisher = {Springer},
}

@article{fletcher2004principal,
  title={Principal geodesic analysis for the study of nonlinear statistics of shape},
  author={Fletcher, P Thomas and Lu, Conglin and Pizer, Stephen M and Joshi, Sarang},
  journal={IEEE transactions on medical imaging},
  volume={23},
  number={8},
  pages={995--1005},
  year={2004},
  publisher={IEEE}
}

@book{durrett2019probability,
  title={Probability: theory and examples},
  author={Durrett, Rick},
  volume={49},
  year={2019},
  publisher={Cambridge university press}
}

@article{jung2012analysis,
  title={Analysis of principal nested spheres},
  author={Jung, Sungkyu and Dryden, Ian L and Marron, James Stephen},
  journal={Biometrika},
  volume={99},
  number={3},
  pages={551--568},
  year={2012},
  publisher={Oxford University Press}
}

@article{korf2001random,
  title={Random lsc functions: an ergodic theorem},
  author={Korf, Lisa A and Wets, Roger J-B},
  journal={Mathematics of Operations Research},
  volume={26},
  number={2},
  pages={421--445},
  year={2001},
  publisher={INFORMS}
}

@article{huckemann2011intrinsic,
  title={Intrinsic inference on the mean geodesic of planar shapes and tree discrimination by leaf growth},
  author={Huckemann, Stephan F},
  journal={Annals of statistics},
  volume={39},
  number={2},
  pages={1098--1124},
  year={2011},
  publisher={Institute of Mathematical Statistics}
}

@inproceedings{ziezold1977expected,
  title={On expected figures and a strong law of large numbers for random elements in quasi-metric spaces},
  author={Ziezold, Herbert},
  booktitle={Transactions of the Seventh Prague Conference on Information Theory, Statistical Decision Functions, Random Processes and of the 1974 European Meeting of Statisticians},
  pages={591--602},
  year={1977},
  organization={Springer}
}

@article{beer1985convergence,
  title={On convergence of closed sets in a metric space and distance functions},
  author={Beer, Gerald},
  journal={Bulletin of the Australian Mathematical Society},
  volume={31},
  number={3},
  pages={421--432},
  year={1985},
  publisher={Cambridge University Press}
}

@article{schotz2022strong,
  title={Strong laws of large numbers for generalizations of {F}r{\'e}chet mean sets},
  author={Sch{\"o}tz, Christof},
  journal={Statistics},
  volume={56},
  number={1},
  pages={34--52},
  year={2022},
  publisher={Taylor \& Francis}
}

@book{beer1993topologies,
  title={Topologies on closed and closed convex sets},
  author={Beer, Gerald},
  volume={268},
  year={1993},
  publisher={Springer Science \& Business Media}
}

@article{afsari2011riemannian,
  title={Riemannian $L^{p}$ center of mass: existence, uniqueness, and convexity},
  author={Afsari, Bijan},
  journal={Proceedings of the American Mathematical Society},
  volume={139},
  number={2},
  pages={655--673},
  year={2011}
}

@article{frechet1948elements,
author = {Fréchet, Maurice},
journal = {Annales de l'institut Henri Poincaré},
language = {fre},
number = {4},
pages = {215-310},
publisher = {INSTITUT HENRI POINCARÉ ET GAUTHIER-VILLARS},
title = {Les éléments aléatoires de nature quelconque dans un espace distancié},
url = {http://eudml.org/doc/79021},
volume = {10},
year = {1948},
}

@incollection{arnaudon2012medians,
  title={Medians and means in {R}iemannian geometry: existence, uniqueness and computation},
  author={Arnaudon, Marc and Barbaresco, Fr{\'e}d{\'e}ric and Yang, Le},
  booktitle={Matrix Information Geometry},
  pages={169--197},
  year={2012},
  publisher={Springer}
}

@article{huckemann2018backward,
author = {Stephan F. Huckemann and Benjamin Eltzner},
title = {{Backward nested descriptors asymptotics with inference on stem cell differentiation}},
volume = {46},
journal = {The Annals of Statistics},
number = {5},
publisher = {Institute of Mathematical Statistics},
pages = {1994 -- 2019},
year = {2018},
doi = {10.1214/17-AOS1609},
URL = {https://doi.org/10.1214/17-AOS1609}
}

@article{jung2025averaging,
  title={Averaging symmetric positive-definite matrices on the space of eigen-decompositions},
  author={Jung, Sungkyu and Rooks, Brian and Groisser, David and Schwartzman, Armin},
  journal={Bernoulli},
  volume={31},
  number={2},
  pages={1552--1578},
  year={2025},
  publisher={Bernoulli Society for Mathematical Statistics and Probability}
}

@article{KIM2024107989,
title = {Principal component analysis for zero-inflated compositional data},
journal = {Computational Statistics \& Data Analysis},
volume = {198},
pages = {107989},
year = {2024},
issn = {0167-9473},
doi = {https://doi.org/10.1016/j.csda.2024.107989},
url = {https://www.sciencedirect.com/science/article/pii/S0167947324000732},
author = {Kipoong Kim and Jaesung Park and Sungkyu Jung}
}

@article{bigot2017geodesic,
author = {J{\'e}r{\'e}mie Bigot and Ra{\'u}l Gouet and Thierry Klein and Alfredo L{\'o}pez},
title = {{Geodesic PCA in the Wasserstein space by convex PCA}},
volume = {53},
journal = {Annales de l'Institut Henri Poincaré, Probabilités et Statistiques},
number = {1},
publisher = {Institut Henri Poincaré},
pages = {1 -- 26},
year = {2017},
doi = {10.1214/15-AIHP706},
URL = {https://doi.org/10.1214/15-AIHP706}
}

@book{howes2012modern,
  title={Modern analysis and topology},
  author={Howes, Norman R},
  year={2012},
  publisher={Springer Science \& Business Media},
pages={38-42}
}

@book{kaufman2009finding,
  title={Finding groups in data: an introduction to cluster analysis},
  author={Kaufman, Leonard and Rousseeuw, Peter J},
  year={2009},
  publisher={John Wiley \& Sons}
}

@article{Huber1964,
author = {Peter J. Huber},
title = {{Robust Estimation of a Location Parameter}},
volume = {35},
journal = {The Annals of Mathematical Statistics},
number = {1},
publisher = {Institute of Mathematical Statistics},
pages = {73 -- 101},
year = {1964},
doi = {10.1214/aoms/1177703732},
URL = {https://doi.org/10.1214/aoms/1177703732}
}

@InProceedings{brunel2023geodesically,
  author       = {Brunel, Victor-Emmanuel},
  booktitle    = {The Thirty Sixth Annual Conference on Learning Theory},
  title        = {Geodesically convex $ M $-estimation in metric spaces},
  year         = {2023},
  organization = {PMLR},
  pages        = {2188--2210},
}

@article{blanchard2025frechet,
  title={Fr{\'e}chet mean set estimation in the Hausdorff metric, via relaxation},
  author={Blanchard, Mo{\"\i}se and Jaffe, Adam Quinn},
  journal={Bernoulli},
  volume={31},
  number={1},
  pages={432--456},
  year={2025},
  publisher={Bernoulli Society for Mathematical Statistics and Probability}
}

@book{srivastava1998course,
  title={A course on Borel sets},
  author={Srivastava, Sashi Mohan},
  year={1998},
  publisher={Springer}
}

@article{evans2024limit,
  title={Limit theorems for Fr{\'e}chet mean sets},
  author={Evans, Steven N and Jaffe, Adam Q},
  journal={Bernoulli},
  volume={30},
  number={1},
  pages={419--447},
  year={2024},
  publisher={Bernoulli Society for Mathematical Statistics and Probability}
}

@book{walters2000introduction,
  title={An introduction to ergodic theory},
  author={Walters, Peter},
  volume={79},
  year={2000},
  publisher={Springer Science \& Business Media}
}
\end{document}